# Computational coupled large-deformation periporomechanics for dynamic failure and fracturing in variably saturated porous media

Shashank Menon, Xiaoyu Song[1]

*Department of Civil and Coastal Engineering, University of Florida, Gainesville, FL-32611.*

**Abstract**

The large-deformation mechanics and multiphysics of continuous or fracturing partially saturated porous media under static and dynamic loads are significant in engineering and science. This article is devoted to a computational coupled large-deformation periporomechanics paradigm assuming passive air pressure for modeling dynamic failure and fracturing in variably saturated porous media. The coupled governing equations for bulk and fracture material points are formulated in the current/deformed configuration through the updated Lagrangian-Eulerian framework. It is hypothesized that the horizon of a mixed material point remains spherical, and its neighbor points are determined in the current configuration. As a significant contribution, the mixed interface/phreatic material points near the phreatic line are explicitly considered for modeling the transition from partial to full saturation (vice versa) through the mixed peridynamic state concept. We have formulated the coupled constitutive correspondence principle and stabilization scheme in the updated Lagrangian-Eulerian framework for bulk and interface points. We numerically implement the coupled large deformation periporomechanics through a fully implicit fractional-step algorithm in time and a hybrid updated Lagrangian-Eulerian meshfree method in space. Numerical examples are presented to validate the implemented stabilized computational coupled large deformation periporomechanics and demonstrate its efficacy and robustness in modeling dynamic failure and fracturing in variably saturated porous media.

## 1. Introduction

The large-deformation mechanics and multiphysics of deformable porous geomaterials (clay and sand) and biomaterials (human tissues and bone) are significant in engineering and science (e.g., geohazards and biomechanics) (e.g., [1–7]). The coupled large deformation and fluid flow/cracking can compromise the integrity of civil infrastructure and could damage human tissues under certain circumstances (e.g., [8–11]). For instance, faulting propagation triggered by earthquake involves large deformation and multiphysics in geomaterials [12–15]. Landslides and landfill slope failure could be triggered by large deformation, multiphysics processes, or cracking in geomaterials [16–18]. Computational coupled poromechanics is an essential tool in studying the mechanics and physics of such continuous or fracturing porous materials under static and dynamic loads (e.g., [2, 19–21]). Coupled unsaturated periporomechanics [22–25] is a strong nonlocal reformulation of classical coupled local poromechanics [26–28] through the peridynamic effective force state [22, 29, 30] and stabilized multiphase correspondence principle [25, 31]. The coupled governing equations, including the motion equation and mass balance equation, are formulated in terms of integral-

[1] Corresponding author

*Email address:* xysong@ufl.edu (Xiaoyu Song )

differential equations (integration in space and differentiation in time) in lieu of partial differential equations [24]. Through the stabilized multiphase correspondence principle, classical constitutive models and physics laws can be incorporated in periporomechanics for modeling the coupled deformation, fracturing, and fluid flow processes in porous media under static and dynamic loads [25]. These salient features of periporomechanics make it a legitimate method for modeling coupled static and dynamic large-deformation mechanics and physics of continuous or fracturing porous media. We note that the coupled periporomechanics has been formulated using the total Lagrangian-Eulerian framework [24]. In [32], the authors formulated an updated-Lagrangian periporomechanics framework for modeling extreme large-deformation in unsaturated porous media under drained conditions (i.e., uncoupled). As a new contribution, in this article, we propose a fully coupled large-deformation periporomechanics paradigm in the updated Lagrangian-Eulerian framework for modeling dynamic failure and fracturing in variably saturated porous media. In this new framework, the phreatic interface/line [28, 33] is explicitly considered through the mixed peridynamic state concept.

To incorporate the classical constitutive models for porous media into this new coupled periporomechanics paradigm, we reformulate the original multiphase constitutive correspondence principle [22] in the updated Lagrangian-Eulerian framework. It has been demonstrated in [25] that the multiphase correspondence principle in the total Lagrangian-Eulerian periporomechanics framework has a zero-energy mode instability issue. In [32] we have also shown that the single-phase correspondence principle for the extreme large-deformation periporomechanics formulated in the updated Lagrangian-Eulerian framework inherits the zero-energy mode instability. The authors adopted the so-called sub-horizon concept to remove the zero-energy mode associated with the correspondence principle. However, this method is computationally demanding because a sub-horizon will need to be defined for each bond in its horizon (see [32] for details). In this study, we first demonstrate that the coupled updated Lagrangian-Eulerian periporomechanics also inherits the zero-energy mode issue (see section 2.3). To reduce the computational cost, we adopt the multiphase stabilization scheme formulated in [25] to resolve the zero-energy mode instability. We refer to the literature for a thorough review of the remedies for circumventing the zero-energy mode associated with the peridynamic correspondence principle (e.g., [25, 31, 32, 34], among others).

The interface between saturated zone and unsaturated zone in variably saturated porous media is called the phreatic interface/line [9, 35]. The porous material across the phreactic line has different mechanical and physical properties because of the variation of degree of saturation across the phreatic line. We note that the phreatic line has not been explicitly considered in the previously formulated satuated/unsaturated periporomechanics models. For a mixed material point near the phreatic line, its neighboring material points can be in either a saturated or unsaturated zone. Therefore, to simulate variably saturated soils with a phreatic line using coupled periporomechanics, the material points near the phreatic line should be treated as composite material points. As a significant contribution and novelty, in this study, the mixed interface/phreatic material points near the phreatic line are explicitly considered for better modeling the transition from partially to fully saturated states (vice versa) of porous media through the mixed peridynamic state concept [29]. Specifically, the peridynamic state (e.g., effective force state and fluid flow state) at the material point across the interface line is decomposed into two states, i.e., saturated and unsaturated states (see Section 2.2 for details). Following the lines in the stabilization scheme for the bulk mixed material point (i.e., material points and their neighboring points in either fully saturated or unsaturated zones), we have formulated the coupled constitutive correspondence principle and its stabilization scheme for the phreatic interface material points in the updated Lagrangian-Eulerian framework.

We numerically implement the coupled large deformation periporomechanics through a fully implicit fractional-step algorithm in time [24] and a hybrid updated Lagrangian-Eulerian meshfree method in space with Open MPI [36] for high-performance computing. We refer to the classical literature for technical discussions on monolithic and fractional-step/staggered algorithms for numerically implementing the coupled poromechanics in time (e.g., [37–41], among others). It is hypothesized that the horizon of a mixed material point remains spherical, and its neighbor points are determined in the current configuration. In line with this hypothesis, in the numerical implementation, the neighboring material list of a material point updated at time step $n$ is used for the computation at time step $n + 1$. The coupled periporomechanics is computationally more demanding than other continuum-based computational methods such as the finite element method (FEM) and XFEM for modeling porous media [1, 2, 19, 20]. We refer to the literature for coupling peridynamics with FEM for modeling porous media (e.g.,[42, 43] and others) and the mixed finite element methods for large deformation in unsaturated porous media (e.g., [2, 21, 44] and others). Numerical examples are presented to validate the implemented stabilized coupled large-deformation periporomechanics and demonstrate its efficacy and robustness in modeling dynamic failure and fracturing in partially saturated porous media.

The original contributions of this article include (i) the mathematical formulation of a fully coupled large-deformation periporomechanics paradigm through the updated Lagrangian-Eulerian framework, (ii) the explicit treatment of interface material points as a composite point through the mixed peridynamic state concept, (iii) the formulation of the coupled constitutive correspondence principle in the updated Lagrangian-Eulerian framework and its stabilization scheme through an energy method, and (iv) the numerical implementation of the computational large-deformation periporomechanics via an implicit-implicit fractional step and mixed meshfree algorithm, and (v) the validation and demonstration of the robustness of the large-deformation periporomechanics paradigm for modeling dynamic failure/fracturing in porous media. For sign convention, the assumption in continuum mechanics is followed, i.e., for the solid skeleton, tensile force/stress is positive, and compression is negative. For fluid pressure, compression is positive, and tension is negative.

## 2. Coupled large-deformation unsaturated periporomechanics

We present the mathematical formulation of the coupled large-deformation periporomechanics in the updated Lagrangian-Eulerian framework for fracturing unsaturated porous media assuming passive air pressure. The phreatic interface points are explicitly considered in this new framework through the mixed peridynamic state concept. In line with the total Lagrangian-Eulerian periporomechanics, it is assumed that a porous material body can be represented by a finite number of mixed material points that are endowed with two kinds of degrees of freedom, i.e., displacement and fluid pressure. In the current/deformed configuration, a mixed material point has poromechanical and physical interactions with all mixed material points within its neighborhood $H$ , i.e., the family. Here $H$ is a spherical domain around a material point $\boldsymbol{x}$ with radius $\delta$, called horizon in periporomechanics, in the current deformation. In the updated Lagrangian-Eulerian formulation, it is hypothesized that the horizon remains the same. In line with this hypothesis, the family $H$ of a mixed material point $\boldsymbol{x}$ in the current configuration is determined by

$$\mathscr{H} := \{\boldsymbol{x}' | \boldsymbol{x}' \in \mathscr{B}, 0 \leq |\underline{\zeta}| \leq \delta\}, \tag{1}$$

where $B$ denotes a porous media body and $\boldsymbol{\zeta} = \boldsymbol{x'} - \boldsymbol{x}$ is the mixed (multiphase) bond between material points $x$ and $x^0$. With this hypothesis, the extreme distortion of the horizon for large deformation of the solid skeleton in the total Lagrangian formulation can be avoided. Indeed, this hypothesis is consistent with the Eulerian formulation of peridynamics for solids in [45]. However, the material point of the solid skeleton is described by its motion, and the fluid phase is described by the relative Eulerian coordinate referring to the solid skeleton. Thus, the mixed material points in the horizon of a material point $\boldsymbol{x}$ evolves with time in the large deformation regime.

Figure 1 presents the schematics of the initial configuration, current/deformed configuration, and future configuration of a porous material body. We note that all variables refer to the current/deformed configuration in the formulation presented in this section instead of the initial/undeformed configuration as in the total Lagrangian-Eulerian periporomechanics. For conciseness of notations, in the current configuration the peridynamic state variable without a prime denotes the variable evaluated at $\boldsymbol{x}$ on the associated bond $\boldsymbol{\zeta} = \boldsymbol{x'} - \boldsymbol{x}$ and the peridynamic state variable with a prime stands for the variable evaluated at $\boldsymbol{x}^0$ on the associated bond $\boldsymbol{\zeta}^0 = \boldsymbol{x} - \boldsymbol{x}^0$, e.g., $T = T[\boldsymbol{x}]h\boldsymbol{x}^0 - \boldsymbol{x}i$ and $T^0 = T[\boldsymbol{x}^0]h\boldsymbol{x} - \boldsymbol{x}^0i$. In what follows, we first present the governing equations and coupled correspondence principle in the updated Lagrangian-Eulerian framework for a mixed material point (and its neighbor material points) (i.e., the bulk material point) in the unsaturated or fully saturated domain. Second, we formulate the governing equations and coupled correspondence principle for the phreatic material point through the mixed peridynamic state concept. Third, we present the stabilization scheme for the phreatic material point since the stabilization scheme for a bulk material point is a special case, followed by a brief discussion of the material models and physical laws in the new framework.

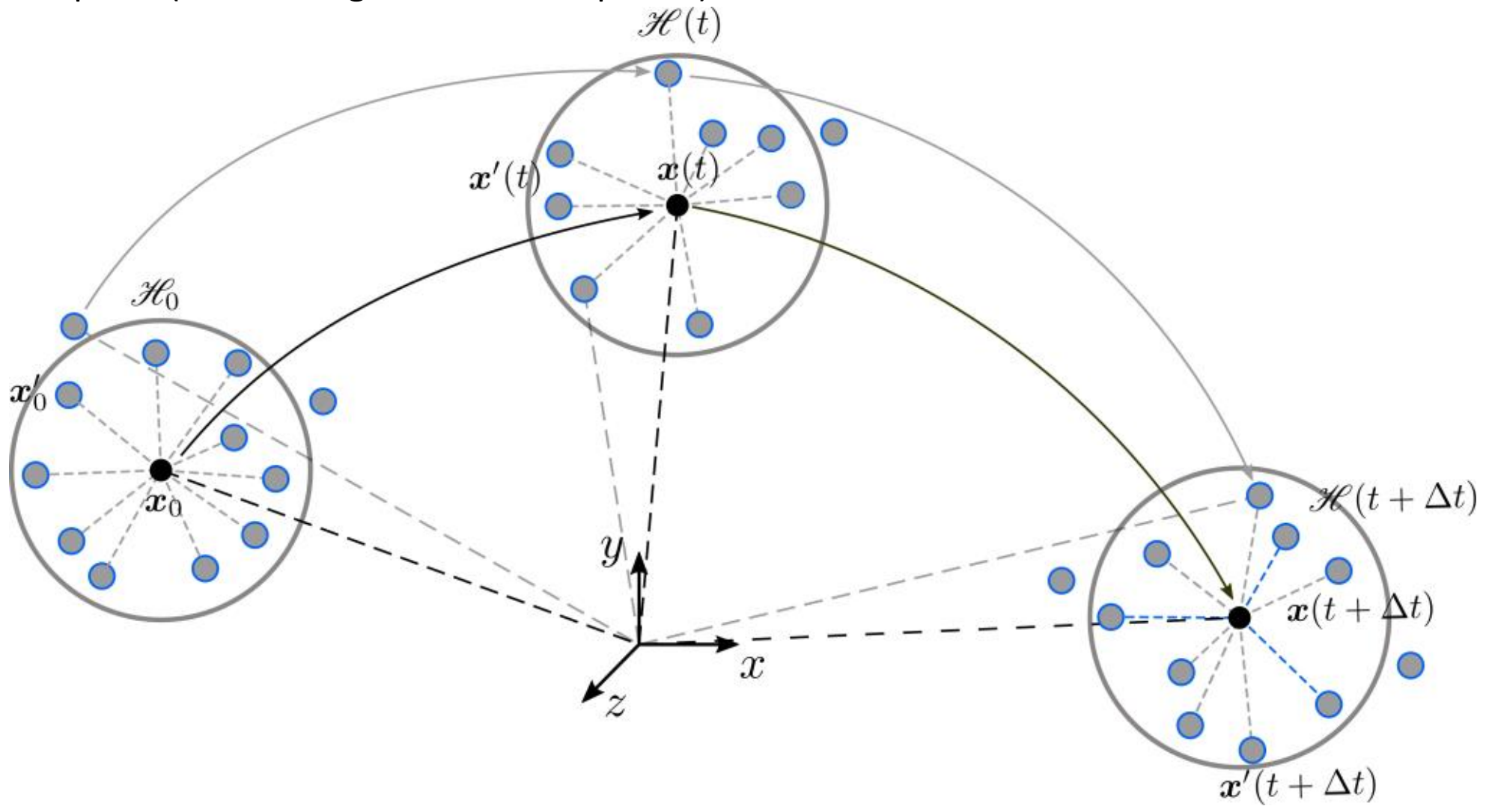

Figure 1: Kinematics of the mixed solid skeleton and pore fluid material points in the updated Lagrangian-Eulerian formulation: initial configuration, current/deformed configuration, and future configuration.

### *2.1. Formulation of a material point in unsaturated or fully saturated zones*

We define the deformation vector state $\underline{Y} = \boldsymbol{y}' - \boldsymbol{y}$ as the mapping of $\boldsymbol{\zeta}$ from the current to the next configuration in which $\boldsymbol{y}'$ and $\boldsymbol{y}$ represent the same material points of the solid skeleton in the next configuration.

In line with the total Lagrangian formulation of periporomechanics [23–25], the equation of motion at point $\boldsymbol{x}$ in the current configuration reads

$$\rho\ddot{\boldsymbol{u}} = \int_{\mathscr{H}} (\underline{\mathscr{T}} - S_r\underline{\mathscr{T}}_w) - (\underline{\mathscr{T}}' - S_r'\underline{\mathscr{T}}_w')\, \mathrm{d}\mathscr{V}' + \rho\boldsymbol{g}, \tag{2}$$

where $\underline{T}$ and $\underline{T}'$ are the effective force vector states at material points $\boldsymbol{x}$ and $\boldsymbol{x}^0$, respectively, $\underline{T}_w$ and $\underline{T}'_w$ are the fluid force vector states at material points $\boldsymbol{x}$ and $\boldsymbol{x}^0$, $\mathrm{d}V^0$ is the volume of the neighboring point in the current domain, and $\ddot{\boldsymbol{u}}$ is the acceleration vector, and $\rho$ is the density of the mixture. Assuming weightless pore air, the density of the mixture $\rho$ reads

$$\rho = \rho_s(1-\phi) + S_r\rho_w\phi, \tag{3}$$

where $\rho_s$ is the intrinsic density of the solid, $\rho_w$ is the intrinsic density of water, $S_r$ is the degree of saturation, and $\varphi$ is the porosity in the current configuration. The balance of mass at $\boldsymbol{x}$ in the current configuration reads

$$\phi\left(\frac{S_r}{K_w} + \frac{\partial S_r}{\partial p_w}\right)\frac{\mathrm{d}p_w}{\mathrm{d}t} + S_r\dot{\mathscr{V}}_s + \int_{\mathscr{H}} \underline{\mathscr{Q}} - \underline{\mathscr{Q}}'\, \mathrm{d}\mathscr{V}' + \mathscr{Q}_s = 0, \tag{4}$$

where $\dot{V}_s$ is the solid volume change rate at $\boldsymbol{x}$, $\underline{Q}$ and $\underline{Q}^0$ are the fluid flow states at $\boldsymbol{x}$ and $\boldsymbol{x}^0$, respectively, $p_w$ is the pore water pressure, $Q_s$ is the source term and $K_w$ is the bulk modulus of water.

For a fracture material point, similar to (2) the equation of motion can be written as

$$\rho \ddot{\mathbf{u}} = \int_{\mathscr{H}} \left( \underline{\overline{\mathscr{T}}} - \underline{\overline{\mathscr{T}}}' \right) \, \mathrm{d}V' - \int_{\mathscr{H}} \left( S_l \underline{\mathscr{T}}_l - S_l' \underline{\mathscr{T}}_l' \right) \, \mathrm{d}V' + \rho \boldsymbol{g}, \tag{5}$$

Where

$$S_l \underline{\mathscr{T}}_l = \begin{cases} S_{r,f} \underline{\mathscr{T}}_f & \text{if } \varphi \,\&\, \varphi' \geq \varphi_{\mathrm{cr}} \\ S_r \underline{\mathscr{T}}_w & \text{otherwises} \end{cases}. \tag{6}$$

The mass balance in the fracture space can be written as

$$\frac{\partial S_{r,f}}{\partial p_w} \frac{\partial p_w}{\partial t} + \int_{\mathscr{H}} \left( \underline{\mathscr{Q}}_f - \underline{\mathscr{Q}}_f' \right) \mathrm{d}\mathscr{V}' - \mathscr{Q}_s = 0, \tag{7}$$

where $\underline{Q}_f$ and $\underline{Q}^0_f$ are fracture flow states at material points $\boldsymbol{x}$ and $\boldsymbol{x}^0$, respectively. In the fully saturated regime, $S_r = 1$ and $\partial S_r/\partial p_w = 0$. Therefore, equations (2),(4), (5), and (7) naturally degenerate into the equations under saturated conditions. In what follows we cast the coupled corresponding principle in the updated-Lagrangian-Eulerian framework.

*2.1.1. Coupled corresponding principle in the updated Lagrangian-Eulerian framework*

The rate of strain energy of the solid skeleton assuming elastic deformation in the updated Lagrangian-Eulerian periporomechanics reads [22]

$$\dot{\mathscr{W}} = \int_{\mathscr{B}} \int_{\mathscr{B}} \underline{\mathscr{T}}_i \, \underline{\dot{\mathscr{Y}}}_i \, \mathrm{d}\mathscr{V}' \mathrm{d}\mathscr{V}. \tag{8}$$

Next we express the rate of strain energy of the solid skeleton in terms of effective stress and the rate of deformation.

The velocity gradient [32] can be written as

$$\boldsymbol{\mathscr{L}} = \left( \int_{\mathscr{H}} \underline{\omega} \, \underline{\dot{\boldsymbol{\mathscr{Y}}}} \otimes \underline{\boldsymbol{\zeta}} \, \mathrm{d}\mathscr{V}' \right) (\boldsymbol{\mathscr{K}})^{-1}, \tag{9}$$

Where K is the spatial shape tensor

$$\boldsymbol{\mathscr{K}} = \int_{\mathscr{H}} \underline{\omega} \, \underline{\boldsymbol{\zeta}} \otimes \underline{\boldsymbol{\zeta}} \, \mathrm{d}\mathscr{V}'. \tag{10}$$

Given (19), the rate of nonlocal deformation reads,

$$\boldsymbol{\mathscr{D}} = \frac{1}{2} [\boldsymbol{\mathscr{L}} + \boldsymbol{\mathscr{L}}^T]. \tag{11}$$

The rate of the strain energy of the solid skeleton assuming elastic deformation can be written as [22, 29],

$$\begin{aligned} \dot{\mathscr{W}} &= \int_{\mathscr{B}} \sigma_{ij} \mathscr{D}_{ij} \, \mathrm{d}\mathscr{V} \\ &= \int_{\mathscr{B}} \sigma_{ij} \left( \int_{\mathscr{H}} \underline{\omega \dot{\mathscr{Y}}}_i \underline{\zeta}_k \, \mathrm{d}\mathscr{V}' \right) \mathscr{K}_{kj}^{-1} \, \mathrm{d}\mathscr{V} \\ &= \int_{\mathscr{B}} \left( \int_{\mathscr{H}} \underline{\omega \dot{\mathscr{Y}}}_i \underline{\zeta}_k \, \mathrm{d}\mathscr{V}' \right) \mathscr{K}_{kj}^{-1} \sigma_{ji} \, \mathrm{d}\mathscr{V} \\ &= \int_{\mathscr{B}} \int_{\mathscr{B}} \underline{\omega \zeta}_k \mathscr{K}_{kj}^{-1} \sigma_{ji} \underline{\dot{\mathscr{Y}}}_i \, \mathrm{d}\mathscr{V}' \, \mathrm{d}\mathscr{V}, \end{aligned} \tag{12}$$

where $i,j,k = 1,2,3$.

Then, it follows from (12) and (8) that the effective force state can be related to the effective Cauchy stress tensor as

$$\underline{\overline{\mathscr{T}}} = \underline{\omega}\underline{\boldsymbol{\zeta}}\boldsymbol{\mathscr{K}}^{-1}\overline{\boldsymbol{\sigma}}, \tag{13}$$

Through the effective force concept, the fluid force state can be written as

$$\underline{\mathscr{T}}_w = \underline{\omega}\underline{\boldsymbol{\zeta}}\boldsymbol{\mathscr{K}}^{-1}(S_r p_w \mathbf{1}). \tag{14}$$

The rate from of energy dissipation due to fluid flow in periporomechanics can be written as [22]

$$\dot{\mathscr{W}}_{\rm d} = \int_{\mathscr{B}} \int_{\mathscr{B}} \underline{\mathscr{Q}}\ \underline{\Phi}\ {\rm d}\mathscr{V}', \tag{15}$$

where $\underline{\Phi} = p^0_w - p_w$ is the pressure potential state. Similar to (9), the spatial gradient of fluid pressure in the current configuration can be written as [23, 32]

$$\mathrm{grad}(p_w)_i = \mathrm{grad}(\widetilde{\Phi})_i = \left(\int_{\mathscr{H}} \underline{\omega}\ \underline{\Phi}\ \underline{\zeta}_j\ {\rm d}\mathscr{V}'\right)\mathscr{K}^{-1}_{ji} \tag{16}$$

The rate form of the energy dissipation due to fluid flow in classical poromechanics [22] reads

$$\begin{aligned}
\dot{\overline{\mathscr{W}}}_d &= \int_{\mathscr{B}} q_i\ \mathrm{grad}(p_w)_i\ {\rm d}\mathscr{V} \\
&= \int_{\mathscr{B}} q_i \left(\int_{\mathscr{H}} \underline{\omega}\ \underline{\Phi}\ \underline{\zeta}_j\ {\rm d}\mathscr{V}'\right)\mathscr{K}^{-1}_{ji}{\rm d}\mathscr{V} \\
&= \int_{\mathscr{B}} q_i \left(\int_{\mathscr{B}} \underline{\omega}\ \underline{\Phi}\ \underline{\zeta}_j\ {\rm d}\mathscr{V}'\right)\mathscr{K}^{-1}_{ji}{\rm d}\mathscr{V} \\
&= \int_{\mathscr{B}} \int_{\mathscr{B}} \underline{\omega} q_i \mathscr{K}^{-1}_{ij} \underline{\zeta}_j \underline{\Phi}\ {\rm d}\mathscr{V}'{\rm d}\mathscr{V},
\end{aligned} \tag{17}$$

Then, from (15) and (17) the fluid flow state can be written as

$$\underline{\mathscr{Q}} = \underline{\omega} q_i \mathscr{K}^{-1}_{ij} \underline{\zeta}_j. \tag{18}$$

Next, we present the governing equations for a mixed material point at the phreatic interface/line.

*2.2. Formulation of a mixed material point at the phreatic interface/line*

In this formulation, the mixed material point near the phreatic line is called the phreatic mixed material points. The phreatic material point could lie in the unsaturated (vadose) zone or the saturated (phreatic) zone. However, a percentage of the material points in its horizon lies in the vadose zone and the remaining material points lie in the phreatic zone. For this reason, it is hypothesized that the deformation state, effective force state, fluid pressure state, and fluid flow state can be decomposed into two parts through the mixed (double) peridynamic state concept, i.e., $\underline{T}^{(1)}$ and $\underline{T}^{(2)}$, $\underline{Y}^{(1)}$ and $\underline{Y}^{(2)}$, $\underline{Q}^{(1)}$ and $\underline{Q}^{(2)}$, and $\underline{\Phi}^{(1)}$ and $\underline{\Phi}^{(2)}$. Figure 2 sketches the periporomechanics representation of the phreatic interface points for the solid skeleton and fluid phases through mixed (double) peridynamic state concept [29].

It follows from this hypothesis and in line with (8) that the strain energy density at a phreatic mixed point reads,

$$\dot{\mathscr{W}} = \int_{\mathscr{B}} \int_{\mathscr{B}} \underline{\overline{\mathscr{T}}}^{(1)}_i \underline{\dot{\mathscr{Y}}}^{(1)}_i + \underline{\overline{\mathscr{T}}}^{(2)}_i \underline{\dot{\mathscr{Y}}}^{(2)}_i\ {\rm d}\mathscr{V}'{\rm d}\mathscr{V}. \tag{19}$$

The two velocity gradient tensors $L^{(1)}$ and $L_{(2)}$ at a phreatic mixed material point can be written as

$$\boldsymbol{\mathscr{L}}^{(1)} = \left( \int_{\mathscr{H}^{(1)}} \underline{\omega}\, \dot{\underline{\boldsymbol{\mathscr{Y}}}} \otimes \underline{\boldsymbol{\zeta}}\ \mathrm{d}\mathscr{V}' \right) \boldsymbol{\mathscr{K}}^{-1}, \tag{20}$$

$$\boldsymbol{\mathscr{L}}^{(2)} = \left( \int_{\mathscr{H}^{(2)}} \underline{\omega}\, \dot{\underline{\boldsymbol{\mathscr{Y}}}} \otimes \underline{\boldsymbol{\zeta}}\ \mathrm{d}\mathscr{V}' \right) \boldsymbol{\mathscr{K}}^{-1}, \tag{21}$$

From (20) and (21) we can obtain the rate form of the two deformation tensor at the interface point, i.e., $D^{(1)}$ and $D^{(2)}$. From a classical constitutive model for the solid skeleton we have two stress tensors, i.e., $\overline{\boldsymbol{\sigma}}^{(1)}$ and $\overline{\boldsymbol{\sigma}}^{(2)}$. Then the corresponding strain energy density of (19) in classical unsaturated poromechanics reads

$$\begin{aligned}
\dot{\mathscr{W}} &= \dot{\mathscr{W}}^{(1)} + \dot{\mathscr{W}}^{(2)} \\
&= \int_{\mathscr{B}} \overline{\sigma}_{ij}^{(1)} \mathscr{D}_{ij}^{(1)}\ \mathrm{d}\mathscr{V} + \int_{\mathscr{B}} \overline{\sigma}_{ij}^{(2)} \mathscr{D}_{ij}^{(2)}\ \mathrm{d}\mathscr{V} \\
&= \int_{\mathscr{B}} \left\{ \overline{\sigma}_{ij}^{(1)} \left( \int_{\mathscr{H}^{(1)}} \underline{\omega}\, \dot{\underline{\mathscr{Y}}}_i \underline{\zeta}_k\ \mathrm{d}\mathscr{V}' \right) \mathscr{K}_{kj}^{-1} \right\} \mathrm{d}\mathscr{V} + \int_{\mathscr{B}} \left\{ \overline{\sigma}_{ij}^{(2)} \left( \int_{\mathscr{H}^{(2)}} \underline{\omega}\, \dot{\underline{\mathscr{Y}}}_i \underline{\zeta}_k\ \mathrm{d}\mathscr{V}' \right) \mathscr{K}_{kj}^{-1} \right\} \mathrm{d}\mathscr{V} \\
&= \int_{\mathscr{B}} \int_{\mathscr{H}^{(1)}} \left\{ \underline{\omega}\, \dot{\underline{\mathscr{Y}}}_i \underline{\zeta}_k \mathscr{K}_{kj}^{-1} \overline{\sigma}_{ji} \right\} \mathrm{d}\mathscr{V}' \mathrm{d}\mathscr{V} + \int_{\mathscr{B}} \int_{\mathscr{H}^{(2)}} \left\{ \underline{\omega}\, \dot{\underline{\mathscr{Y}}}_i \underline{\zeta}_k \mathscr{K}_{kj}^{-1} \overline{\sigma}_{ji} \right\} \mathrm{d}\mathscr{V}' \mathrm{d}\mathscr{V} \\
&= \int_{\mathscr{B}} \int_{\mathscr{H}} \left\{ \underline{\omega}\, \dot{\underline{\mathscr{Y}}}_i \underline{\zeta}_k \mathscr{K}_{kj}^{-1} \overline{\sigma}_{ji} \right\}^{(1)} \mathrm{d}\mathscr{V}' \mathrm{d}\mathscr{V} + \int_{\mathscr{B}} \int_{\mathscr{H}} \left\{ \underline{\omega}\, \dot{\underline{\mathscr{Y}}}_i \underline{\zeta}_k \mathscr{K}_{kj}^{-1} \overline{\sigma}_{ji} \right\}^{(2)} \mathrm{d}\mathscr{V}' \mathrm{d}\mathscr{V} \\
&= \int_{\mathscr{B}} \int_{\mathscr{B}} \left( \underline{\omega}\, \underline{\zeta}_k \mathscr{K}_{kj}^{-1} \overline{\sigma}_{ji} \right)^{(1)} \dot{\underline{\mathscr{Y}}}_i^{(1)} + \left( \underline{\omega}\, \underline{\zeta}_k \mathscr{K}_{kj}^{-1} \overline{\sigma}_{ji} \right)^{(2)} \dot{\underline{\mathscr{Y}}}_i^{(2)}\ \mathrm{d}\mathscr{V}' \mathrm{d}\mathscr{V},
\end{aligned} \tag{27}$$

where $\dot{W}^{(1)}$ and $\dot{W}^{(2)}$ are strain energy densities associated with the neighboring material points under saturated and unsaturated conditions, respectively.

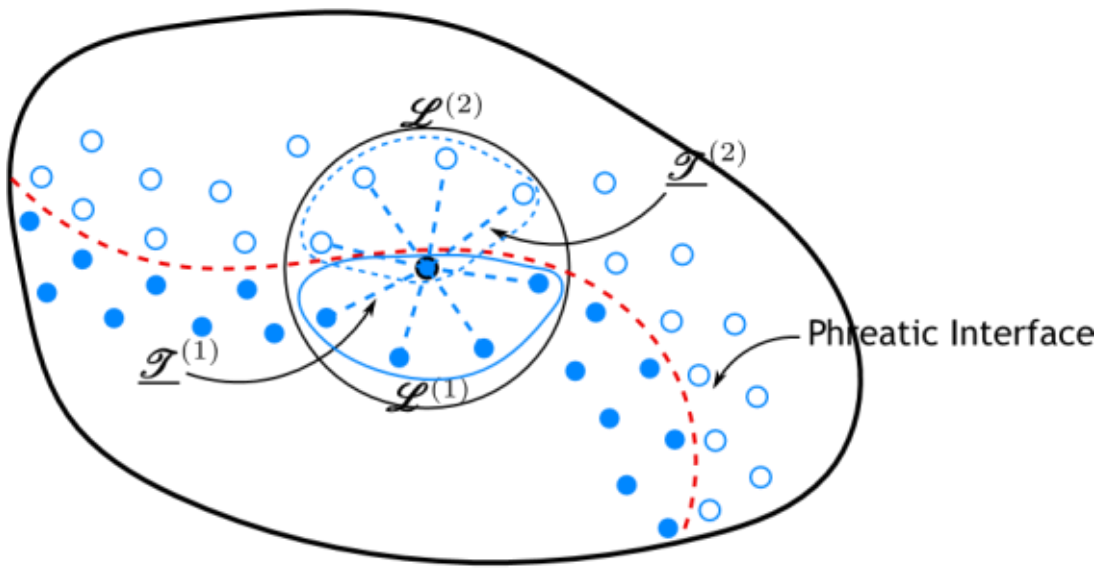

(a)

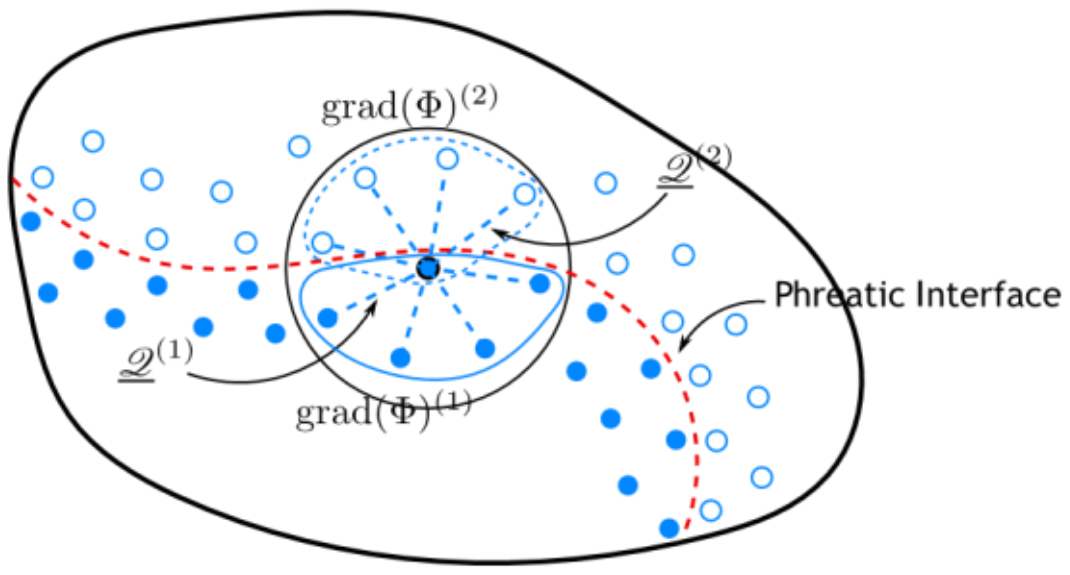

(b)

Figure 2: Schematics of the phreatic interface material points in periporomechanics for (a) the solid skeleton phase and (b) the pore fluid phase through mixed (double) peridynamic state concept.

From (19) and (22), it follows that the effective force states at a phreatic material point can be written as

$$\underline{\overline{\mathscr{T}}}^{(1)} = \left\{ \underline{\omega} \zeta \boldsymbol{\mathscr{K}}^{-1} (\overline{\boldsymbol{\sigma}}) \right\}^{(1)} \tag{23}$$

$$\underline{\overline{\mathscr{T}}}^{(2)} = \left\{ \underline{\omega} \zeta \boldsymbol{\mathscr{K}}^{-1} (\overline{\boldsymbol{\sigma}}) \right\}^{(2)} \tag{24}$$

Through the effective force state concept, the fluid force states at the phreactic material point read

$$\underline{\mathscr{T}}_w^{(1)} = \left\{ \underline{\omega} \zeta \boldsymbol{\mathscr{K}}^{-1} (S_r p_w \mathbf{1}) \right\}^{(1)} \tag{25}$$

$$\underline{\mathscr{T}}_w^{(2)} = \left\{ \underline{\omega} \zeta \boldsymbol{\mathscr{K}}^{-1} (S_r p_w \mathbf{1}) \right\}^{(2)} \tag{26}$$

Then, substituting (23), (24), (25) and (26) into (2) we have the equation of motion at the phreatic material point in terms of the total force states as

$$\rho \ddot{\boldsymbol{u}} = \int_{\mathscr{H}} \left( \underline{\mathscr{T}}^{(1)} + \underline{\mathscr{T}}^{(2)} \right) - \left( \underline{\mathscr{T}}'^{(1)} + \underline{\mathscr{T}}'^{(2)} \right) \, \mathrm{d}\mathscr{V}' + \rho \boldsymbol{g}. \tag{27}$$

For the fluid phase, the energy dissipation rate at a phreatic material point through fluid flow can be written as

$$\dot{\mathscr{W}}_d = \int_{\mathscr{B}} \int_{\mathscr{B}} \left( \underline{\mathscr{Q}}^{(1)} \, \underline{\Phi}^{(1)} + \underline{\mathscr{Q}}^{(2)} \, \underline{\Phi}^{(2)} \right) \, \mathrm{d}\mathscr{V}'. \tag{28}$$

At a phreatic material point, the two fluid pressure gradient vectors can be written as

$$\widetilde{\mathrm{grad}(\Phi)}^{(1)} = \left( \int_{\mathscr{H}^{(1)}} \underline{\omega} \, \underline{\Phi} \, \underline{\zeta} \, \mathrm{d}\mathscr{V}' \right) \boldsymbol{\mathscr{K}}^{-1}, \tag{29}$$

$$\widetilde{\mathrm{grad}(\Phi)}^{(2)} = \left( \int_{\mathscr{H}^{(2)}} \underline{\omega} \, \underline{\Phi} \, \underline{\zeta} \, \mathrm{d}\mathscr{V}' \right) \boldsymbol{\mathscr{K}}^{-1}. \tag{30}$$

Given (29) and (30), flux vectors $\boldsymbol{q}^{(1)}$ and $\boldsymbol{q}^{(2)}$ can be obtained through the generalized Darcy's law for unsaturated porous media. Then, the energy dissipation rate due to fluid flow at a phreatic material point in classical poromechanics reads

$$\begin{aligned}
\dot{\overline{\mathscr{W}}}_d &= \dot{\overline{\mathscr{W}}}_d^{(1)} + \dot{\overline{\mathscr{W}}}_d^{(2)} \\
&= \int_{\mathscr{B}} q_i^{(1)} \, \mathrm{grad}(p_w)_i^{(1)} \, \mathrm{d}\mathscr{V} + \int_{\mathscr{B}} q_i^{(2)} \, \mathrm{grad}(p_w)_i^{(2)} \, \mathrm{d}\mathscr{V} \\
&= \int_{\mathscr{B}} \left\{ q_i^{(1)} \left( \int_{\mathscr{H}^{(1)}} \underline{\omega} \, \underline{\Phi} \, \underline{\zeta}_j \, \mathrm{d}\mathscr{V}' \right) \mathscr{K}_{ji}^{-1} \right\} \mathrm{d}\mathscr{V} + \int_{\mathscr{B}} \left\{ q_i^{(2)} \left( \int_{\mathscr{H}^{(2)}} \underline{\omega} \, \underline{\Phi} \, \underline{\zeta}_j \, \mathrm{d}\mathscr{V}' \right) \mathscr{K}_{ji}^{-1} \right\} \mathrm{d}\mathscr{V} \\
&= \int_{\mathscr{B}} \int_{\mathscr{H}} \left\{ \underline{\omega} \, \underline{\Phi} \, \underline{\zeta}_j \mathscr{K}_{ji}^{-1} q_i \right\}^{(1)} \mathrm{d}\mathscr{V}' \mathrm{d}\mathscr{V} + \int_{\mathscr{B}} \int_{\mathscr{H}} \left\{ \underline{\omega} \, \underline{\Phi} \, \underline{\zeta}_j \mathscr{K}_{ji}^{-1} q_i \right\}^{(2)} \mathrm{d}\mathscr{V}' \mathrm{d}\mathscr{V} \\
&= \int_{\mathscr{B}} \int_{\mathscr{B}} \left( \underline{\omega} \, \underline{\zeta}_j \mathscr{K}_{ji}^{-1} q_i \right)^{(1)} \underline{\Phi}^{(1)} + \left( \underline{\omega} \, \underline{\zeta}_j \mathscr{K}_{ji}^{-1} q_i \right)^{(2)} \underline{\Phi}^{(2)} \, \mathrm{d}\mathscr{V}' \mathrm{d}\mathscr{V}.
\end{aligned} \tag{36}$$

It follows from (31) and (28) that the two fluid flow states at a phreatic material points read

$$\underline{\mathscr{Q}}^{(1)} = \left\{ \underline{\omega} \boldsymbol{q} \boldsymbol{\mathscr{K}}^{-1} \underline{\boldsymbol{\zeta}} \right\}^{(1)}, \tag{32}$$

$$\underline{\mathscr{Q}}^{(2)} = \left\{ \underline{\omega} \boldsymbol{q} \boldsymbol{\mathscr{K}}^{-1} \underline{\boldsymbol{\zeta}} \right\}^{(2)}. \tag{33}$$

Given (32) and (33) the balance of mass at a phreatic material point can be written as

$$\phi \left( \frac{S_r}{K_w} + \frac{\partial S_r}{\partial p_w} \right) \frac{\mathrm{d}p_w}{\mathrm{d}t} + S_r \sum_{i=1}^{2} \underline{\mathscr{L}}^{(i)} : \mathbf{1} + \int_{\mathscr{H}} \left[ \left( \underline{\mathscr{Q}}^{(1)} + \underline{\mathscr{Q}}^{(2)} \right) - \left( \underline{\mathscr{Q}}'^{(1)} + \underline{\mathscr{Q}}^{(2)} \right) \right] \mathrm{d}\mathscr{V}' + \mathscr{Q}_s = 0. \tag{34}$$

In this study we assume that the fracturing is not across the phreatic interface/line. However, the formulation for fracturing material points across the phreatic line can be achieved following the lines in this section and in [23].

*2.3. Coupled stabilization scheme in the updated Lagrangian-Eulerian framework*

In this section, we derive an two-phase stabilized scheme for the two-phase correspondence principle for the phreatic material point derived in the previous section. Note that the stabilization scheme of the mixed material point completely in either saturated or unsaturated domain is a special case of the phreatic material point. For a phreatic material point, the two spatial deformation gradient vectors are written as

$$\boldsymbol{\mathscr{F}}^{(1)} = \int_{\mathscr{H}^{(1)}} \underline{\omega}\ \underline{\boldsymbol{\mathscr{Y}}} \otimes \underline{\boldsymbol{\zeta}}\ \mathrm{d}\mathscr{V}' \boldsymbol{\mathscr{K}}^{-1}, \tag{35}$$

$$\boldsymbol{\mathscr{F}}^{(2)} = \int_{\mathscr{H}^{(2)}} \underline{\omega}\ \underline{\boldsymbol{\mathscr{Y}}} \otimes \underline{\boldsymbol{\zeta}}\ \mathrm{d}\mathscr{V}' \boldsymbol{\mathscr{K}}^{-1}. \tag{36}$$

The two residual deformation vector states that are written as

$$\underline{\boldsymbol{\mathscr{R}}}^{s,(1)} = \underline{\boldsymbol{\mathscr{Y}}}^{(1)} - \boldsymbol{\mathscr{F}}^{(1)} \underline{\boldsymbol{\zeta}}^{(1)}, \tag{37}$$

$$\underline{\boldsymbol{\mathscr{R}}}^{s,(2)} = \underline{\boldsymbol{\mathscr{Y}}}^{(2)} - \boldsymbol{\mathscr{F}}^{(2)} \underline{\boldsymbol{\zeta}}^{(2)}. \tag{38}$$

Substituting (37) and (38) into (35) and (36) generates two null tensors implying that the nonuniform deformation state is not considered by the correspondence principle derived in the last section [32]. The missing non-uniform deformation state is the origin of the zero-energy mode associated with the standard correspondence principle.

The energy method [25] is adopted for its computational efficiency. The effective force state with stabilization can be written as

$$\overline{\underline{\boldsymbol{\mathscr{T}}}} = \overline{\underline{\boldsymbol{\mathscr{T}}}}^{(1)} + \overline{\underline{\boldsymbol{\mathscr{T}}}}^{(2)} + \underline{\boldsymbol{\mathscr{T}}}_s^{(1)} + \underline{\boldsymbol{\mathscr{T}}}_s^{(2)}, \tag{39}$$

where $\underline{T}^{(1)}{}_s$ and $\underline{T}^{(2)}{}_s$ are the stabilization terms corresponding to $\overline{\underline{T}}^{(1)}$ and $\underline{T}^{(2)}$, respectively. Next, we determine the two stabilization terms from the energy method. Following [25, 31] we propose that the stabilization force states take the form

$$\underline{\boldsymbol{\mathscr{T}}}_s^{(1)} = \beta^{(1)} \underline{\boldsymbol{\mathscr{R}}}^{s,(1)}, \tag{40}$$

$$\underline{\boldsymbol{\mathscr{T}}}_s^{(2)} = \beta^{(2)} \underline{\boldsymbol{\mathscr{R}}}^{s,(2)}, \tag{41}$$

where $\beta^{(1)}$ and $\beta^{(2)}$ are weighting factors. The strain energy associated with the non-uniform deformation is written as

$$\begin{aligned} \mathscr{W}_s(\underline{\boldsymbol{\mathscr{Y}}}) &= \mathscr{W}_s(\underline{\boldsymbol{\mathscr{Y}}}^{(1)}) + \mathscr{W}_s(\underline{\boldsymbol{\mathscr{Y}}}^{(2)}), \\ &= \frac{1}{2}(\underline{\beta \boldsymbol{\mathscr{R}}})^{s,(1)} \bullet \underline{\boldsymbol{\mathscr{R}}}^{s,(1)} + \frac{1}{2}(\underline{\beta \boldsymbol{\mathscr{R}}})^{s,(2)} \bullet \underline{\boldsymbol{\mathscr{R}}}^{s,(2)}, \end{aligned} \tag{42}$$

where • is the state dot product operator. $W_s(\underline{Y}^{(1)})$ in (42) can be expressed as

$$\begin{aligned} \mathscr{W}_s^{(1)} &= \int_{\mathscr{H}^{(1)}} \underline{\beta}^{(1)} \underline{\mathscr{R}}_i^{s,(1)} (\underline{\mathscr{Y}}_i - \mathscr{F}_{ij}\zeta_j)^{(1)}\, \mathrm{d}\mathscr{V}' \\ &= \int_{\mathscr{H}^{(1)}} \{\underline{\beta \mathscr{R}_i^{s}} \underline{\mathscr{Y}}_i\}^{(1)}\ \mathrm{d}\mathscr{V}' - \int_{\mathscr{H}^{(1)}} \{\underline{\beta \mathscr{R}_i^{s}} \mathscr{F}_{ij}\zeta_j\}^{(1)}\ \mathrm{d}\mathscr{V}' \\ &= \int_{\mathscr{H}^{(1)}} \{\underline{\beta \mathscr{R}_i^{s}} \underline{\mathscr{Y}}_i\}^{(1)}\ \mathrm{d}\mathscr{V}' - \left(\int_{\mathscr{H}^{(1)}} \underline{\beta \mathscr{R}_i^{s}} \zeta_j \mathrm{d}\mathscr{V}'\right)^{(1)} \left(\int_{\mathscr{H}^{(1)}} \underline{\omega \mathscr{Y}}_i \zeta_l\ \mathrm{d}\mathscr{V}'\right)^{(1)} \mathscr{K}_{lj}^{-1} \\ &= \int_{\mathscr{H}^{(1)}} \left[\underline{\beta \mathscr{R}_i^{s}} - \left(\int_{\mathscr{H}^{(1)}} \underline{\beta \mathscr{R}_i^{s}} \zeta_j\ \mathrm{d}\mathscr{V}'\right) \underline{\omega} \zeta_l \mathscr{K}_{lj}^{-1}\right]^{(1)} \underline{\mathscr{Y}}_i^{(1)} \mathrm{d}\mathscr{V}'. \end{aligned} \tag{43}$$

Similarly, we have

$$\mathscr{W}_s^{(2)} = \int_{\mathscr{H}^{(2)}} \left[\underline{\beta \mathscr{R}_i^{s}} - \left(\int_{\mathscr{H}^{(1)}} \underline{\beta \mathscr{R}_i^{s}} \zeta_j\ \mathrm{d}\mathscr{V}'\right) \underline{\omega} \zeta_l \mathscr{K}_{lj}^{-1}\right]^{(2)} \underline{\mathscr{Y}}_i^{(2)} \mathrm{d}\mathscr{V}'. \tag{44}$$

In line with the formulation of the standard correspondence principle the stabilized force state terms at a phreatic material point can be determined as

$$\underline{\mathscr{T}}^{(1)}_{s,i} = \left\{ \underline{\beta \mathscr{R}}^s_i - \left( \int_{\mathscr{H}} \underline{\beta \mathscr{R}}^s_i \zeta_j \, \mathrm{d}\mathscr{V}' \right) \underline{\omega} \zeta_l \mathscr{K}^{-1}_{lj} \right\}^{(1)}, \tag{45}$$

$$\underline{\mathscr{T}}^{(2)}_{s,i} = \left\{ \underline{\beta \mathscr{R}}^s_i - \left( \int_{\mathscr{H}} \underline{\beta \mathscr{R}}^s_i \zeta_j \, \mathrm{d}\mathscr{V}' \right) \underline{\omega} \zeta_l \mathscr{K}^{-1}_{lj} \right\}^{(2)}. \tag{46}$$

In this study to simplify (45) and (46) it is assumed that $\beta^{(1)}$ and $\beta^{(2)}$ are constants as

$$\beta^{(1)} = \frac{GC^{(1)}}{\omega_0^{(1)}} \underline{\omega}, \tag{47}$$

$$\beta^{(2)} = \frac{GC^{(2)}}{\omega_0^{(2)}} \underline{\omega}, \tag{48}$$

where $G$ is the stabilization parameter [25, 31], $C$ is the micro-elastic modulus [25], and

$$\omega_0^{(1)} = \int_{\mathscr{H}^{(1)}} \underline{\omega} \, \mathrm{d}\mathscr{V}',$$

$$\omega_0^{(2)} = \int_{\mathscr{H}^{(2)}} \underline{\omega} \, \mathrm{d}\mathscr{V}'. \tag{49}$$

It follows then that (45) and (46) can be written as

$$\underline{\boldsymbol{\mathscr{T}}}^{s,(1)} = \frac{GC^{(1)}}{\omega_0^{(1)}} \underline{\omega \boldsymbol{\mathscr{R}}}^{s,(1)}, \tag{50}$$

$$\underline{\boldsymbol{\mathscr{T}}}^{s,(2)} = \frac{GC^{(2)}}{\omega_0^{(2)}} \underline{\omega \boldsymbol{\mathscr{R}}}^{s,(2)}. \tag{51}$$

Then the effective force state at a phreatic material point can be written as

$$\overline{\underline{\boldsymbol{\mathscr{T}}}} = \underline{\omega} \left\{ \zeta \boldsymbol{\mathscr{K}}^{-1} \overline{\boldsymbol{\sigma}} + \frac{GC}{\omega_0} \underline{\boldsymbol{\mathscr{R}}}^s \right\}^{(1)} + \underline{\omega} \left\{ \zeta \boldsymbol{\mathscr{K}}^{-1} \overline{\boldsymbol{\sigma}} + \frac{GC}{\omega_0} \underline{\boldsymbol{\mathscr{R}}}^s \right\}^{(2)}. \tag{52}$$

The parameters $C^{(1)}$ and $C^{(2)}$ can be determined by the strain energy method assuming an isotropic elastic deformation of the skeleton as follows. For simplicity, a microelastic peridynamic model [46] is adopted to determine the elastic energy in the solid skeleton.

Let $\boldsymbol{u}$ and $\boldsymbol{u}^0$ are the displacement vectors at material points $\boldsymbol{x}$ and $\boldsymbol{x}^0$ respectively referring to the current configuration. The rate form of the strain energy density at point $\boldsymbol{x}$ due to deformation $\boldsymbol{\eta} = \boldsymbol{u}^0 - \boldsymbol{u}$ can be written as

$$\begin{aligned} \dot{\mathscr{W}} &= \frac{1}{2} \int_{\mathscr{H}^{(1)}} w(\eta, \zeta) \, \mathrm{d}\mathscr{V}' + \frac{1}{2} \int_{\mathscr{H}^{(2)}} w(\eta, \zeta) \, \mathrm{d}\mathscr{V}' \\ &= \frac{1}{2} \varphi^{(1)} \int_{\mathscr{H}} w(\eta, \zeta) \, \mathrm{d}\mathscr{V}' + \frac{1}{2} \varphi^{(2)} \int_{\mathscr{H}} w(\eta, \zeta) \, \mathrm{d}\mathscr{V}', \end{aligned} \tag{53}$$

where $\eta = ||\boldsymbol{\eta}||$, $\zeta = ||\boldsymbol{\zeta}||$, and $\phi^{(1)}$ and $\phi^{(2)}$ are defined as follows.

$$\varphi^{(1)} = \frac{\int_{\mathscr{H}^{(1)}} \mathrm{d}\mathscr{V}'}{\int_{\mathscr{H}} \mathrm{d}\mathscr{V}'}, \tag{54}$$

$$\varphi^{(2)} = \frac{\int_{\mathscr{H}^{(2)}} \mathrm{d}\mathscr{V}'}{\int_{\mathscr{H}} \mathrm{d}\mathscr{V}'}. \tag{55}$$

Consider an isotropic deformation of the solid skeleton the stretch of a bond reads

$$\eta = \mathsf{C}_1\zeta, \tag{56}$$

where $C_1$ is a material constant. The micro-potential [25] can be written as

$$w^{(1)} = C^{(1)}\mathsf{C}_1^2\zeta/2, \tag{57}$$

$$w^{(2)} = C^{(2)}\mathsf{C}_1^2\zeta/2. \tag{58}$$

Substituting (57) and (58) into (53) gives

$$\begin{aligned}\dot{\mathscr{W}} &= \frac{1}{2}\varphi^{(1)}\int_{\mathscr{H}}\frac{1}{2}\left(C^{(1)}\dot{\mathsf{C}}_1^2\zeta\right)\,\mathrm{d}\mathscr{V}' + \frac{1}{2}\varphi^{(2)}\int_{\mathscr{H}}\frac{1}{2}\left(C^{(2)}\dot{\mathsf{C}}_1^2\zeta\right)\,\mathrm{d}\mathscr{V}' \\ &= \frac{1}{2}\varphi^{(1)}\int_0^{\delta}\frac{1}{2}\left(C^{(1)}\dot{\mathsf{C}}_1^2\zeta\right)(4\pi\zeta^2)\,\mathrm{d}\zeta + \frac{1}{2}\varphi^{(2)}\int_0^{\delta}\frac{1}{2}\left(C^{(2)}\dot{\mathsf{C}}_1^2\zeta\right)(4\pi\zeta^2)\,\mathrm{d}\zeta \\ &= \varphi^{(1)}\frac{\pi C^{(1)}\delta^4}{4}\dot{\mathsf{C}}_1^2 + \varphi^{(2)}\frac{\pi C^{(2)}\delta^4}{4}\dot{\mathsf{C}}_1^2.\end{aligned} \tag{59}$$

From the classical elastic theory, the incremental elastic strain energy density of the solid skeleton at a phreatic material point under isotropic deformation reads

$$\begin{aligned}\dot{\widetilde{\mathscr{W}}} &= \left\{\frac{1}{2}\overline{\boldsymbol{\sigma}}:\dot{\boldsymbol{\epsilon}}\right\}^{(1)} + \left\{\frac{1}{2}\overline{\boldsymbol{\sigma}}:\dot{\boldsymbol{\epsilon}}\right\}^{(2)} \\ &= \left\{\frac{1}{2}(K^e\epsilon_v)(\dot{\epsilon}_v)\right\}^{(1)} + \left\{\frac{1}{2}(K^e\epsilon_v)(\dot{\epsilon}_v)\right\}^{(2)} \\ &= \left\{\frac{1}{2}K^e(3\dot{\mathsf{C}}_1)^2\right\}^{(1)} + \left\{\frac{1}{2}K^e(3\dot{\mathsf{C}}_1)^2\right\}^{(2)} \\ &= \left\{\frac{9K^e}{2}\dot{\mathsf{C}}_1^2\right\}^{(1)} + \left\{\frac{1}{2}K^e(3\dot{\mathsf{C}}_1)^2\right\}^{(2)},\end{aligned} \tag{60}$$

where *K* is the elastic bulk modulus and $_v$ is the elastic volumetric strain under isotropic deformation. From (60) and (59) *C* under three-dimensional condition can be written as

$$C^{(1)} = \varphi^{(1)}\frac{18K^e}{\pi\delta^4}, \tag{61}$$

$$C^{(2)} = \varphi^{(2)}\frac{18K^e}{\pi\delta^4}, \tag{62}$$

Next, we first demonstrate the fluid flow state in the updated Lagrangian-Eulerian formulation inherits the instability as in the total Lagrangian-Eulerian formulation and then present an stabilization scheme through the energy method for the fluid phase. The two residual flow states at a phreatic material points read

$$\mathscr{R}^{w,(1)} = \underline{\Phi}^{(1)} - \mathrm{grad}\widetilde{(\Phi)}^{(1)}\underline{\boldsymbol{\zeta}}^{(1)}, \tag{63}$$

$$\mathscr{R}^{w,(2)} = \underline{\Phi}^{(2)} - \mathrm{grad}\widetilde{(\Phi)}^{(2)}\underline{\boldsymbol{\zeta}}^{(2)}. \tag{64}$$

Substituting (63) or (64) into the spatial nonlocal pressure gradient (16) generates

$$
\begin{aligned}
\operatorname{grad}\left(\underline{\mathscr{R}}^{w}\right) &= \left(\int_{\mathscr{H}} \underline{\omega}\ \underline{\Phi}\ \underline{\boldsymbol{\zeta}}\ \mathrm{d}\mathscr{V}'\right) \boldsymbol{\mathscr{K}}^{-1} \\
&= \int_{\mathscr{H}} \underline{\omega}\ \underline{\mathscr{R}}^{w}\ \underline{\boldsymbol{\zeta}}\ \mathrm{d}\mathscr{V}' \left(\int_{\mathscr{H}} \underline{\omega}\ \underline{\boldsymbol{\zeta}} \otimes \underline{\boldsymbol{\zeta}}\ \mathrm{d}\mathscr{V}'\right)^{-1} \\
&= \int_{\mathscr{H}} \underline{\omega}\ (\underline{\Phi} - \operatorname{grad}(\widetilde{\Phi})\underline{\boldsymbol{\zeta}})\underline{\boldsymbol{\zeta}}\ \mathrm{d}\mathscr{V}' \left(\int_{\mathscr{H}} \underline{\omega}\ \underline{\boldsymbol{\zeta}} \otimes \underline{\boldsymbol{\zeta}}\ \mathrm{d}\mathscr{V}'\right)^{-1} \\
&= \left(\int_{\mathscr{H}} \underline{\omega}\ (\underline{\Phi}\underline{\boldsymbol{\zeta}})\ \mathrm{d}\mathscr{V}' - \operatorname{grad}(\widetilde{\Phi}) \int_{\mathscr{H}} \underline{\omega}\ (\underline{\boldsymbol{\zeta}} \otimes \underline{\boldsymbol{\zeta}})\ \mathrm{d}\mathscr{V}'\right) \left(\int_{\mathscr{H}} \underline{\omega}\ \underline{\boldsymbol{\zeta}} \otimes \underline{\boldsymbol{\zeta}}\ \mathrm{d}\mathscr{V}'\right)^{-1} \\
&= \operatorname{grad}(\widetilde{\Phi}) - \operatorname{grad}(\widetilde{\Phi})\boldsymbol{\mathscr{K}}\boldsymbol{\mathscr{K}}^{-1} = 0.
\end{aligned} \tag{65}
$$

It is implied by (65) the correspondence principle for fluid flow in the updated Lagrangian-Eulerian formulation will generate zero-energy mode instability under non-uniform fluid flow conditions. We define the stabilization terms for the flow states at a phreatic material point as

$$
\underline{\mathscr{Q}}^{s} = \underline{\lambda}^{(1)}\underline{\mathscr{R}}^{w,(1)} + \underline{\lambda}^{(2)}\underline{\mathscr{R}}^{w,(2)}, \tag{66}
$$

where $\underline{\lambda}^{(1)}$ and $\underline{\lambda}^{(2)}$ are positive parameters. The fluid flow dissipation energy through periporomechanics in the updated-Lagrangian formulation reads

$$
\mathscr{W}_{d}^{s} = \frac{1}{2}\left(\underline{\lambda\mathscr{R}}^{w}\right)^{w,(1)} \bullet \underline{\mathscr{R}}^{w,(1)} + \frac{1}{2}\left(\underline{\lambda\mathscr{R}}^{w}\right)^{(2)} \bullet \underline{\mathscr{R}}^{w,(2)}. \tag{67}
$$

It follows from (63) and (67) that

$$
\begin{aligned}
\mathscr{W}_{d}^{s,(1)} &= \int_{\mathscr{H}^{(1)}} \underline{\lambda}^{(1)}\underline{\mathscr{R}}^{w,(1)} \left(\underline{\Phi} - \operatorname{grad}(\widetilde{\Phi})_{i}\underline{\zeta}_{i}\right)^{(1)} \mathrm{d}\mathscr{V}' \\
&= \int_{\mathscr{H}^{(1)}} \{\underline{\lambda\mathscr{R}}^{w}\underline{\Phi}\}^{(1)} \mathrm{d}\mathscr{V}' - \int_{\mathscr{H}^{(1)}} \left\{\underline{\lambda\mathscr{R}}^{w}\operatorname{grad}(\widetilde{\Phi})_{i}\underline{\zeta}_{i}\right\}^{(1)} \mathrm{d}\mathscr{V}' \\
&= \int_{\mathscr{H}^{(1)}} \{\underline{\lambda\mathscr{R}}^{w}\underline{\Phi}\}^{(1)} \mathrm{d}\mathscr{V}' - \left(\int_{\mathscr{H}^{(1)}} \underline{\lambda\mathscr{R}}\underline{\zeta}_{i}\mathrm{d}\mathscr{V}'\right)^{(1)} \left(\int_{\mathscr{H}^{(1)}} \underline{\omega\Phi}\underline{\zeta}_{j}\mathrm{d}\mathscr{V}'\right)^{(1)} \mathscr{K}_{ji}^{-1} \\
&= \int_{\mathscr{H}^{(1)}} \left[\underline{\lambda\mathscr{R}}^{w} - \left(\int_{\mathscr{H}^{(1)}} \underline{\lambda\mathscr{R}}^{w}\underline{\zeta}_{i}\mathrm{d}\mathscr{V}'\right)\underline{\omega}\underline{\zeta}_{j}\mathscr{K}_{ji}^{-1}\right]^{(1)} \underline{\Phi}^{(1)}\ \mathrm{d}\mathscr{V}'.
\end{aligned} \tag{68}
$$

Similarly, we have

$$
\mathscr{W}_{d}^{s,(2)} = \int_{\mathscr{H}^{(1)}} \left[\underline{\lambda\mathscr{R}}^{w} - \left(\int_{\mathscr{H}^{(2)}} \underline{\lambda\mathscr{R}}^{w}\underline{\zeta}_{i}\mathrm{d}\mathscr{V}'\right)\underline{\omega}\underline{\zeta}_{j}\mathscr{K}_{ji}^{-1}\right]^{(2)} \underline{\Phi}^{(2)}\ \mathrm{d}\mathscr{V}'. \tag{69}
$$

From (15), (68) (69) the stabilized terms for the fluid flow state accounting for the nonuniform fluid potential state can be written as

$$
\underline{\mathscr{Q}}^{s,(1)} = \left\{\underline{\lambda\mathscr{R}}^{w} - \left(\int_{\mathscr{H}} \underline{\lambda\mathscr{R}}^{w}\zeta_{i}\mathrm{d}\mathscr{V}'\right)\underline{\omega}\zeta_{j}\mathscr{K}_{ji}^{-1}\right\}^{(1)}, \tag{70}
$$

$$
\underline{\mathscr{Q}}^{s,(2)} = \left\{\underline{\lambda\mathscr{R}}^{w} - \left(\int_{\mathscr{H}} \underline{\lambda\mathscr{R}}^{w}\zeta_{i}\mathscr{K}_{ji}^{-1}\mathrm{d}\mathscr{V}'\right)\underline{\omega}\zeta_{j}\right\}^{(2)}. \tag{71}
$$

For simplicity, we assume that $\lambda$ takes the form

$$\underline{\lambda}^{(1)} = \frac{GK_p^{(1)}}{\omega_0^{(1)}}\underline{\omega}, \tag{72}$$

$$\underline{\lambda}^{(2)} = \frac{GK_p^{(2)}}{\omega_0^{(2)}}\underline{\omega}, \tag{73}$$

Where $K_p^1$ and $K_p^2$ are the hydraulic micro-conductivities. Then we obtain

$$\underline{\mathscr{Q}}^{s,(1)} = \frac{GK_p^{(1)}}{\omega_0^{(1)}}\underline{\omega \mathscr{R}}^{w,(1)},$$

$$\underline{\mathscr{Q}}^{s,(2)} = \frac{GK_p^{(2)}}{\omega_0^{(2)}}\underline{\omega \mathscr{R}}^{w,(2)}. \tag{74}$$

Then, the stabilized fluid flow state at $\boldsymbol{x}$ is written as

$$\underline{\mathscr{Q}} = \underline{\omega}\left(\boldsymbol{q}\boldsymbol{\mathscr{K}}^{-1}\underline{\zeta} + \frac{GK_p}{\omega_0}\underline{\mathscr{R}}^w\right)^{(1)} + \underline{\omega}\left(\boldsymbol{q}\boldsymbol{\mathscr{K}}^{-1}\underline{\zeta} + \frac{GK_p}{\omega_0}\underline{\mathscr{R}}^w\right)^{(2)} \tag{75}$$

Next, we derive the material variables through an energy method assuming an isotropic fluid flow. The dissipation potential energy at a phreatic point reads,

$$\begin{aligned}\mathscr{W}_d &= \frac{1}{2}\int_{\mathscr{H}^{(1)}} w^f(\boldsymbol{x}',\boldsymbol{x})\ \mathrm{d}\mathscr{V}' + \frac{1}{2}\int_{\mathscr{H}^{(2)}} w^f(\boldsymbol{x}',\boldsymbol{x})\ \mathrm{d}\mathscr{V}' \\ &= \frac{1}{2}\varphi^{(1)}\int_{\mathscr{H}} w^f(\boldsymbol{x}',\boldsymbol{x})\ \mathrm{d}\mathscr{V}' + \frac{1}{2}\varphi^{(2)}\int_{\mathscr{H}} w^f(\boldsymbol{x}',\boldsymbol{x})\ \mathrm{d}\mathscr{V}'.\end{aligned} \tag{76}$$

The peridynamic hydraulic micro-conductivity can be related to the classical hydraulic conductivity by equating the peridynamic fluid dissipation potential to the classical fluid dissipation potential at a phreatic point. For simplicity, we assume a linear pressure field in a body, $p_w = C_2(\mathbf{1}\cdot\boldsymbol{x})$ for a three-dimensional case. The fluid flow micro potential [25] can be written as

$$w^f = K_p\frac{1}{2}\mathscr{C}_2^2\zeta. \tag{77}$$

Then we have

$$\begin{aligned}\mathscr{W}_d &= \frac{1}{2}\varphi^{(1)}\int_0^\delta\left(\frac{1}{2}K_p^{(1)}\mathscr{C}_2^2\zeta\right)(4\pi\zeta^2)\mathrm{d}\zeta + \frac{1}{2}\varphi^{(2)}\int_0^\delta\left(\frac{1}{2}K_p^{(2)}\mathscr{C}_2^2\zeta\right)(4\pi\zeta^2)\mathrm{d}\zeta \\ &= \varphi^{(1)}\frac{\pi K_p^{(1)}\delta^4}{4}\mathscr{C}_2^2 + \varphi^{(2)}\frac{\pi K_p^{(2)}\delta^4}{4}\mathscr{C}_2^2.\end{aligned} \tag{78}$$

Assuming an isotropic fluid flow, the classical fluid dissipation energy in classical poromechanics reads

$$\begin{aligned}\overline{\mathscr{W}}_d &= \left\{\frac{1}{2}(\boldsymbol{\nabla}p_w)\frac{k_r k_w}{\mu}\mathbf{1}(\boldsymbol{\nabla}p_w)\right\}^{(1)} + \left\{\frac{1}{2}(\boldsymbol{\nabla}p_w)\frac{k_r k_w}{\mu}\mathbf{1}(\boldsymbol{\nabla}p_w)\right\}^{(2)} \\ &= \left\{\frac{3}{2}\frac{k_r k_w}{\mu}\mathscr{C}_2^2\right\}^{(1)} + \left\{\frac{3}{2}\frac{k_r k_w}{\mu}\mathscr{C}_2^2\right\}^{(2)}.\end{aligned} \tag{79}$$

where $k_w$ is the intrinsic permeability of saturated porous media and $k_r$ is the relative permeability for the partially saturated case.

It follows from (78) and (79) that the hydraulic micro-conductivities at a phreatic material point under three-dimensional condition are

$$K_p^{(1)} = \varphi^{(1)} \frac{6k_r k_w}{\mu\pi\delta^4}, \tag{80}$$

$$K_p^{(2)} = \varphi^{(2)} \frac{6k_r k_w}{\mu\pi\delta^4}. \tag{81}$$

*2.4. Constitutive models and physical laws*

We summarize the constitutive models and physical laws. The effective stress tensor can be determined by an isotropic linear elastic model (or a plastic or visco-plastic constitutive model [32]) as

$$\dot{\hat{\boldsymbol{\sigma}}} = \mathfrak{C} : \mathscr{D}^e, \tag{82}$$

where C is the forth-order elastic modulus tensor [4], $\dot{\hat{\boldsymbol{\sigma}}}$ is the rotated Cauchy stress and $D^e$ is the elastic deformation rate tensor. Given the velocity gradient tensor the rate of unrotated deformation rate tensor reads

$$\mathscr{D} = \mathscr{R}^T \left( \frac{1}{2}[\mathscr{L} + \mathscr{L}^T] \right) \mathscr{R}, \tag{83}$$

where *R* is the rotation tensor from the left polar decomposition of *F*.

$$\mathscr{F} = \mathscr{V}\mathscr{R}.$$

Once the effective Cauchy stress tensor is updated it is rotated back to the deformed configuration as follows

$$\dot{\boldsymbol{\sigma}} = \mathscr{R}^T \dot{\hat{\boldsymbol{\sigma}}} \mathscr{R} \tag{85}$$

Similarly, the fluid flux $\boldsymbol{q}$ can be determined by the generalized Darcy's law for unsaturated fluid flow as

$$\boldsymbol{q} = -\frac{k^r \boldsymbol{k}_w}{\mu_w} \operatorname{grad}(\widetilde{\Phi}), \tag{86}$$

where $k^r$ is the relative permeability, grad$(\Phi)_g$ is the nonlocal fluid pressure gradient referring to the current configuration, $\mu$ is the water viscosity and $\boldsymbol{k}_w$ is the intrinsic permeability tensor.

For large deformation applications there exists pure rotations that change the reference frame with no deformation. Therefore the permeability tensor needs to modified as

$$\boldsymbol{k}_w = \mathscr{R}^T \hat{\boldsymbol{k}}_w \mathscr{R}. \tag{87}$$

It follows from (87) that (86) [2] can be written as

$$\boldsymbol{q} = -\mathscr{R} \left[ \frac{k^r \mathscr{R}^T \boldsymbol{k}_w \mathscr{R}}{\mu_w} \operatorname{grad}(\widetilde{\Phi}) \right]. \tag{88}$$

The fluid flow state in the fracture space reads

$$\underline{Q}_f = \frac{3}{m_v} \omega \rho_w \boldsymbol{q}_f \underline{\zeta}. \tag{89}$$

where $\boldsymbol{q}_f$ is the fluid flow vector in fracture space. Through Darcy's law for unsaturated fluid flow the fracture fluid flow vector $\boldsymbol{q}_f$ can be written as

$$\boldsymbol{q}_f = -\frac{k_f^r k_f}{\mu_w} \widetilde{\boldsymbol{\nabla}\Phi}_f, \tag{90}$$

where $k_f^r$ is the relative permeability, $k_f$ is the intrinsic permeability of fracture space, and $\nabla_g \Phi_f$ is the nonlocal fracture fluid pressure gradient determined as

$$\widetilde{\nabla\Phi}_f = \frac{3}{m_v} \int_{\mathscr{H}} \underline{\omega}\ \underline{\Phi}_f \underline{\zeta}\ \mathrm{d}\mathscr{V}', \tag{91}$$

where

$$\underline{\Phi}_f = p_f' - p_f, \tag{92}$$

And

$$m_v = \int_{\mathscr{H}} \underline{\omega}\ \underline{\zeta}\ \underline{\zeta}\ \mathrm{d}\mathscr{V}'. \tag{93}$$

The porosity in (3) can be written as [21]

$$\phi = 1 - \frac{(1-\phi_t)\mathscr{J}}{\mathscr{J}}. \tag{94}$$

where $\varphi_t$ is the current configuration porosity and

$$J = \det(F) \tag{95}$$

The degree of saturation $S_r$ can be determined from the soil-water retention curve (e.g., [47–51]) that depends on the volume strain of the solid skeleton (e.g., porosity). In this study, we adopt the one in [21, 52, 53] which reads

$$S_r(\mathscr{J}, \phi, p_w) = \left\{1 + \left[-a_1 \left(\frac{\mathscr{J}}{1-\phi} - 1\right)^{a_2} \mathscr{J} p_w\right]^n\right\}^{(n-1)/n}, \tag{96}$$

where $a_1$, $a_2$, and $n$ are all material parameters. Given $S_r$ the relative permeability $k_r$ can be written as

$$k^r = S_r^{1/2}\left[1 - (1 - S_r^{1/m})^m\right]^2, \tag{97}$$

where $m = (n-1)/n$.

The bond-breakage criteria is based on the deformation energy stored in a poromechanical bond. With effective force state the energy density in an intact poromechanical bond $\boldsymbol{\zeta}$ reads

$$\underline{\varpi} = \int_t^{t+\Delta t} \left(\underline{\mathscr{T}} - \underline{\mathscr{T}}'\right) \dot{\boldsymbol{\eta}}\ \mathrm{d}t = \int_t^{t+\Delta t} \left[(\underline{\mathscr{T}} + S_r \underline{\mathscr{T}}_w) - (\underline{\mathscr{T}}' + S_r' \underline{\mathscr{T}}_w')\right] \dot{\boldsymbol{\eta}}\ \mathrm{d}t, \tag{98}$$

where $\dot{\boldsymbol{\eta}} = \dot{\boldsymbol{u}}^0 - \dot{\boldsymbol{u}}$ is the relative displacement vector. The bond breakage is realized through the influence function at the material point level for both the solid and fluid phases. The influence function will be replaced by a new influence function $\%\underline{\omega}$, where $\%$ is defined as

$$\underline{\varrho} = \begin{cases} 0 & , \text{ if } \underline{\varpi} \geq \varpi_{cr} \\ 1 & , \text{ otherwise} \end{cases} \tag{99}$$

In periporomechanics, the failure at solid skeleton material points is tracked through a scalar damage variable $\phi$ [24]. This damage variable is defined as the fraction of broken solid bonds at a material point in its horizon

$$\varphi = 1 - \frac{\int_{\mathscr{H}} \underline{\varrho}\underline{\omega} \, \mathrm{d}\mathscr{V}'}{\omega_0}, \tag{100}$$

where $\phi \in [0,1]$ and $\omega_0$ is defined as

$$\omega_0 = \int_{\mathscr{H}} \underline{\omega} \, \mathrm{d}\mathscr{V}'. \tag{101}$$

## 3. Numerical implementation

In this section, we implement the computational large-deformation periporomechanics model through an implicit-implicit fractional step algorithm. With an undrained operator split, the strongly coupled problem is decomposed into an undrained deformation stage and a partially saturated fluid flow stage. The list of neighboring material points of a mixed material point is updated in the deformed configuration at each time step through a dedicated search algorithm. The evolution of phreatic material points across the phreatic line is tracked at each time step through another dedicated search algorithm. In this study, we focus on deformation-driven processes, and loading conditions are applied to the solid skeleton. Thus, we first solve the deformation problem and then the fluid flow process. We note that in other applications, such as hydraulic fracturing, wherein fluid pressure drives deformation, it would be preferred to solve the fluid flow first.

### *3.1. Time discretization of the deformation stage*

In the deformation stage, the equation of motion is solved under undrained condition through an implicit Newton's method [24, 39].

Given $\boldsymbol{u}_n$, $\dot{\boldsymbol{u}}_n$, $\ddot{\boldsymbol{u}}_n$ at time step $n$. Let $\Delta\ddot{\boldsymbol{u}}^{k+1}_{n+1} = \ddot{\boldsymbol{u}}^{k+1}_{n+1} - \ddot{\boldsymbol{u}}_n$ and $k$ is the iteration counter. Following Newmark's method [1], the displacement, velocity, and accelerations of the solid skeleton at time step $n + 1$ can be written as

$$\ddot{\boldsymbol{u}}^{k+1}_{n+1} = \ddot{\boldsymbol{u}}^{k}_{n+1} + \delta\Delta\ddot{\boldsymbol{u}}^{k+1}_{n+1}, \tag{102}$$

$$\dot{\boldsymbol{u}}^{k+1}_{n+1} = \dot{\boldsymbol{u}}_n + \Delta t\ddot{\boldsymbol{u}}_n + \beta_2\Delta t\Delta\ddot{\boldsymbol{u}}^{k+1}_{n+1}, \tag{103}$$

$$\boldsymbol{u}^{k+1}_{n+1} = \boldsymbol{u}_n + \Delta t\dot{\boldsymbol{u}}_n + \frac{(\Delta t)^2}{2}\ddot{\boldsymbol{u}}_n + \beta_1\frac{(\Delta t)^2}{2}\Delta\ddot{\boldsymbol{u}}^{k+1}_{n+1}, \tag{104}$$

where $\beta_1$ and $\beta_2$ are numerical integration parameters. For unconditional stability [1] $\beta_1 \geq \beta_2 \geq 0.5$.

The solution procedure for the solid deformation stage with phreatic material points is outlined in Figure 3. A detailed description of its numerical implementation is provided in Algorithm 1. At $t_{n+1}$, the residual of the motion equation for a mixed point in the bulk space reads

$$\boldsymbol{r}^{u,k+1} = \rho^{k+1}\ddot{\boldsymbol{u}}^{k+1} - \int_{\mathscr{H}_n} \underline{\varrho}_n \left(\underline{\mathscr{T}} - \underline{\mathscr{T}}'\right)^{k+1} \mathrm{d}\mathscr{V}'_n + \int_{\mathscr{H}_n} \left(S_r\underline{\mathscr{T}}_w - S'_r\underline{\mathscr{T}}'_w\right)^{k+1} \mathrm{d}\mathscr{V}'_n - \rho^{k+1}\boldsymbol{g}. \tag{105}$$

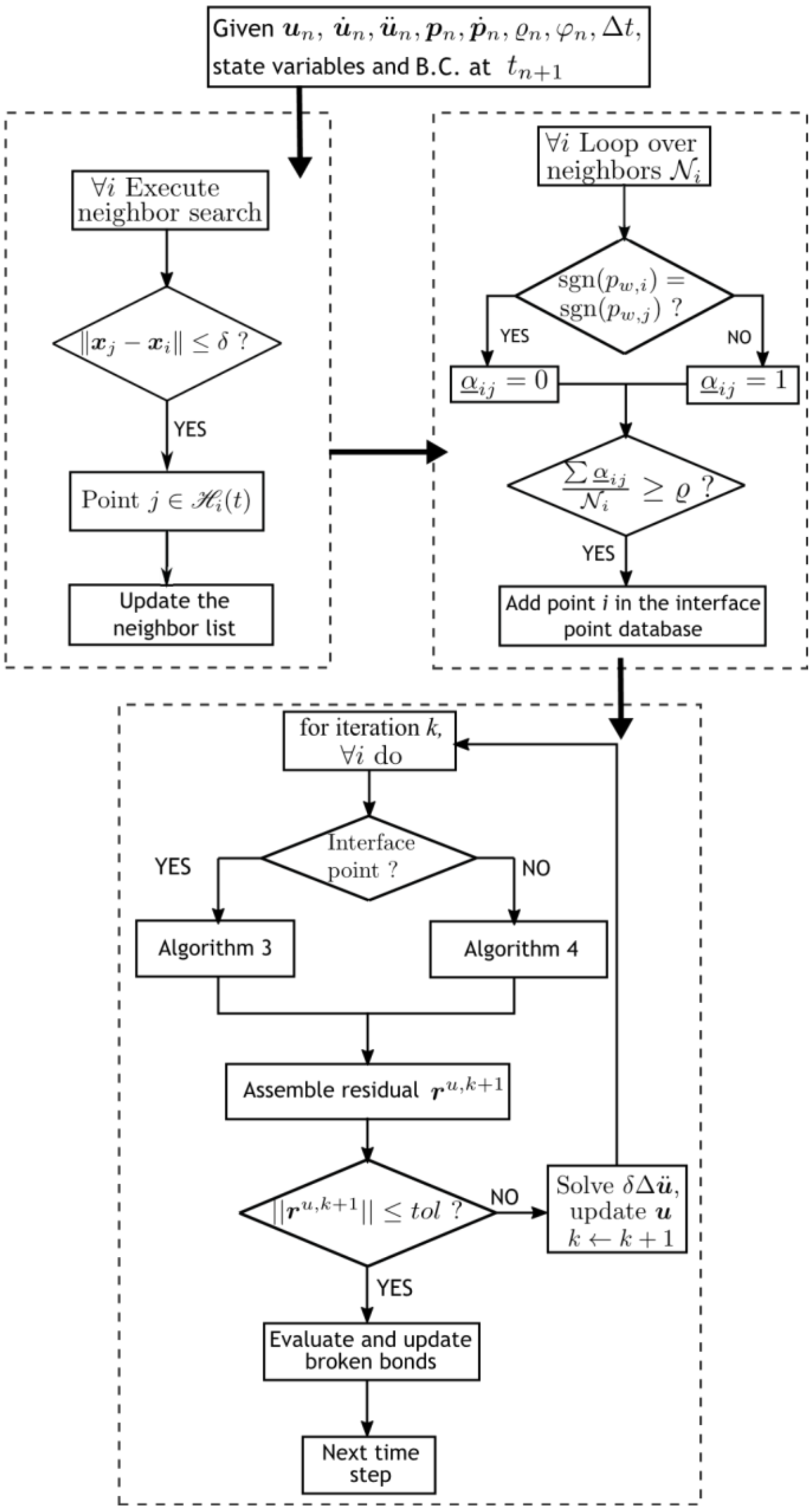

Figure 3: Flowchart of Newton's method for solving the motion equation in the deformation stage of the fractional step algorithm for coupled updated Lagrangian-Eulerian periporomechanics.

**Algorithm 1** Fractional-step algorithm for the solid deformation stage in fully coupled updated Lagrangian-Eulerian periporomechanics

**procedure** GIVEN $\boldsymbol{u}_n, \dot{\boldsymbol{u}}_n, \ddot{\boldsymbol{u}}_n, \boldsymbol{p}_{w,n}, \dot{\boldsymbol{p}}_{w,n}, \boldsymbol{p}_{f,n}, \dot{\boldsymbol{p}}_{f,n}, t_n$ AND $\Delta t$, SOLVE $\ddot{\boldsymbol{u}}_{n+1}$ AND COMPUTE $\overline{\mathscr{V}}_{n+1}$
  $t_{n+1} = t_n + \Delta t$
  **for** all points $i$ **do**
    **for** all points $j$ **do**
      Execute neighbor search $\forall j \in \mathscr{B}$
      **if** $||\boldsymbol{x}_{(j)} - \boldsymbol{x}_{(i)}|| \leq \delta$ **then**
        Add $j \in \mathscr{H}_i$
      **end if**
    **end for**
  **end for**
  Compute the velocity predictor $\tilde{\dot{\boldsymbol{u}}}_{n+1}$
  Apply boundary conditions
  Compute the displacement predictor $\widetilde{\boldsymbol{u}}_{n+1}$
  Compute the predictors $\widetilde{\boldsymbol{p}}_{w,n+1}$ and $\widetilde{\boldsymbol{p}}_{f,n+1}$
  **for** all points $i$ **do**
    **for** all neighbors $j$ **do**
      **if** $sgn(p_{w,i}) = sgn(p_{w,j})$ **then**
        set $\underline{\alpha}_{ij} = 0$
      **else if** $sgn(p_{w,i})! = sgn(p_{w,j})$ **then**
        set $\underline{\alpha}_{ij} = 1$
      **end if**
    **end for**
    **if** $\sum \underline{\alpha}_{ij} / \mathscr{N}_i > \varsigma$ **then**
      SET $i$ as INTERFACE_POINT
    **end if**
  **end for**
  Compute the effective force and fluid force via Algorithm 2
  Compute balance of momentum residual $\mathscr{R}^{u,0}_{n+1}$
  Set k = 0,
  **while** $||\boldsymbol{r}^{u,k}_{n+1}|| > \mathrm{tol}_1$ **do**
    Construct tangent operator for the momentum balance $\mathscr{A}^{u,k} = \left[\partial \boldsymbol{r}^{u,k} / \partial \Delta \ddot{\boldsymbol{u}}^k\right]_{n+1}$
    Solve the linear system for $\delta \Delta \ddot{\boldsymbol{u}}^{k+1}_{n+1}$
    Update $\ddot{\boldsymbol{u}}^{k+1}_{n+1}$, $\boldsymbol{x}^{k+1}_{n+1}$, $\widetilde{\boldsymbol{p}}^{k+1}_{w,n+1}$ and $\widetilde{\boldsymbol{p}}^{k+1}_{f,n+1}$
    Check the residual $\mathscr{R}^{u,k+1}_{n+1}$
    Set $k \leftarrow k+1$
  **end while**
  Update $\boldsymbol{p}_{w,n} \leftarrow \widetilde{\boldsymbol{p}}_{w,n+1}$, $\boldsymbol{p}_{f,n} \leftarrow \widetilde{\boldsymbol{p}}_{f,n+1}$
  Evaluate the bond failure criterion and update broken bonds
**end procedure**

Similarly, the residual of the motion equation for a phreatic material point reads

$$\boldsymbol{r}^{u,k+1} = \rho^{k+1} \ddot{\boldsymbol{u}}^{k+1} - \int_{\mathscr{H}_n} \underline{\varrho}_n \left\{ \left( \underline{\boldsymbol{\mathscr{T}}}^{(1)} + \underline{\boldsymbol{\mathscr{T}}}^{(2)} \right) - \left( \underline{\boldsymbol{\mathscr{T}}}'^{(1)} + \underline{\boldsymbol{\mathscr{T}}}'^{(2)} \right) \right\}^{k+1} \mathrm{d}\mathscr{V}'_n + \int_{\mathscr{H}_n} \left\{ S_r \left( \underline{\boldsymbol{\mathscr{T}}}^{(1)}_w + \underline{\boldsymbol{\mathscr{T}}}^{(2)}_w \right) - S'_r \left( \underline{\boldsymbol{\mathscr{T}}}'^{(1)}_w + \underline{\boldsymbol{\mathscr{T}}}'^{(2)}_w \right) \right\}^{k+1} \mathrm{d}\mathscr{V}'_n - \rho^{k+1} \boldsymbol{g}. \tag{106}$$

The terms of the fluid pressure in (105) are determined from the explicit fluid pressure predictor $\tilde{p}_w^{k+1}$ computed under the undrained condition. Given $p_{w,n}$ and $\dot{p}_{w,n}$ at time step $n$, let $\Delta\dot{\tilde{p}}_{w,n+1} = \dot{\tilde{p}}_{w,n+1} - \dot{p}_{w,n}$. From Newmark's method [1, 25] $\tilde{p}_{w,n+1}$ can be written as

$$\tilde{p}_w^{k+1} = p_{w,n} + \Delta t \dot{p}_{w,n} + \beta_3 \Delta t \Delta \dot{\tilde{p}}^{k+1}, \tag{107}$$

Through the undrained operator split, the rate form of $\tilde{p}_{w,n+1}$ at a material point in the bulk Reads

$$\dot{\tilde{p}}_w^{k+1} = -\left(\phi \frac{S_r}{K_w} + \phi \frac{\partial S_r}{\partial p_w}\right)_n^{-1} \left[S_{r,n}(\boldsymbol{\mathscr{L}}^{k+1} : \mathbf{1}) + \int_{\mathscr{H}_n} \underline{\varrho}_n \left(\underline{\mathscr{Q}}_n - \underline{\mathscr{Q}}'_n\right) \, \mathrm{d}\mathscr{V}'_n\right]. \tag{108}$$

Similarly, the rate of the predicted fluid pressure at a phreatic material point can be written as

$$\begin{aligned} \dot{\tilde{p}}_w^{k+1} = &-\left(\phi \frac{S_r}{K_w} + \phi \frac{\partial S_r}{\partial p_w}\right)_n^{-1} \left[S_{r,n}\left(\boldsymbol{\mathscr{L}}^{(1),k+1} + \boldsymbol{\mathscr{L}}^{(2),k+1}\right) : \mathbf{1}\right. \\ &\left. + \int_{\mathscr{H}_n} \underline{\varrho}_n \left\{\left(\underline{\mathscr{Q}}^{(1)} + \underline{\mathscr{Q}}^{(2)}\right)_n - \left(\underline{\mathscr{Q}}'^{(1)} + \underline{\mathscr{Q}}'^{(2)}\right)_n\right\} \, \mathrm{d}\mathscr{V}'_n\right]. \end{aligned} \tag{109}$$

Algorithms 2 and 3 summarize the algorithms for computing the force states in the residual the motion equation at a bulk material point and a phreatic material point, respectively.

---

**Algorithm 2** Construct $\underline{\boldsymbol{\mathscr{T}}}_{n+1}$ and $\underline{\mathscr{T}}_{w,n+1}$ for interface points

---

**for** all interface points $i$ **do**
  **for** all neighbors $j$ **do**
    Compute the contribution to $\boldsymbol{\mathscr{K}}_{(i)}$
    **if** $\underline{\alpha}_{ij} = 0$ **then**
      Compute the contribution to $\boldsymbol{\mathscr{F}}_i^{k,(1)}$, $\boldsymbol{\mathscr{L}}_i^{k,(1)}$ and $\boldsymbol{\mathscr{D}}_i^{k,(1)}$
    **end if**
    **if** $\underline{\alpha}_{ij} = 1$ **then**
      Compute the contribution to $\boldsymbol{\mathscr{F}}_i^{k,(2)}$, $\boldsymbol{\mathscr{L}}_i^{k,(2)}$ and $\boldsymbol{\mathscr{D}}_i^{k,(2)}$
    **end if**
  **end for**
**end for**
**for** all interface points $i$ **do**
  Update the porosity $\phi_i^k$
  Compute the volume coupling term $\dot{\mathscr{V}}_i^k$
  Compute the rate pressure predictor $\dot{\tilde{p}}_i^k$
  Compute $\Delta\boldsymbol{\varepsilon}_i^{k,(1)}$ and $\overline{\boldsymbol{\sigma}}_i^{k,(1)}$
  Compute $\Delta\boldsymbol{\varepsilon}_i^{k,(2)}$ and $\overline{\boldsymbol{\sigma}}_i^{k,(2)}$
**end for**
**for** all interface points $i$ **do**
  **for** all neighbors $j$ **do**
    **if** $\underline{\alpha}_{ij} = 0$ **then**
      Compute $\underline{\overline{\boldsymbol{\mathscr{T}}}}_{ij}^{k,(1)}$ and $\underline{\mathscr{T}}_{w,ij}^{k,(1)}$
    **end if**
    **if** $\underline{\alpha}_{ij} = 1$ **then**
      Compute $\underline{\overline{\boldsymbol{\mathscr{T}}}}_{ij}^{k,(2)}$ and $\underline{\mathscr{T}}_{w,ij}^{k,(2)}$
    **end if**
  **end for**
**end for**

---

Substituting (102), (103), and (104) into (105) and using (107), $\Delta\ddot{\boldsymbol{u}}_{n+1}$ can be solved by Newton's method as follows.

$$\mathscr{R}^{u,k+1} = \mathscr{R}^{u,k} + \mathscr{A}^{u,k}\delta\Delta\ddot{\boldsymbol{u}}^{k+1} \approx \mathbf{0}, \tag{110}$$

where $\mathbf{R}^u$ is the global residual vector of the motion equations and $\mathbf{A}^u$ is the global tangent operator of the motion equation

$$\mathscr{A}^{u,k} = \frac{\partial \mathscr{R}^{u,k}}{\partial \Delta\ddot{\boldsymbol{u}}^k}. \tag{111}$$

Solving (111) we obtain

$$\delta\Delta\ddot{\boldsymbol{u}}^{k+1} = -(\mathscr{A}^k)^{-1}\mathscr{R}^{u,k}. \tag{112}$$

Finally, we have,

$$\Delta\ddot{\boldsymbol{u}}^{k+1} = \Delta\ddot{\boldsymbol{u}}^k + \delta\Delta\ddot{\boldsymbol{u}}^{k+1}. \tag{113}$$

*3.2. Time discretization of the unsaturated fluid flow stage*

Given $p_{w,n}$, $\dot{p}_{w,n}$, and $\boldsymbol{u}_n$, the unsaturated fluid flow stage solves $\dot{p}_{w,n+1}$ in the updated configuration of the solid skeleton configuration using an implicit Newton's method at time step $n+1$. Let $\Delta\dot{p}^k_w{}^{+1} = \dot{p}^k_w{}^{+1} - \dot{p}_{w,n}$. From Newmark's method [1] we have

$$\dot{p}_w^{k+1} = \dot{p}_w^k + \delta\Delta\dot{p}_w^{k+1}, \tag{114}$$

$$p_w^{k+1} = p_{w,n} + \Delta t\dot{p}_{w,n} + \beta_3\Delta t\Delta\dot{p}_w^{k+1}, \tag{115}$$

The solution procedure for the unsaturated flow stage with phreatic material points is outlined in Figure 4. A detailed description of its numerical implementation is provided in Algorithm 4. At $t_{n+1}$, the residual of the mass balance equation for a bulk material point is written as

$$r^{p,k+1} = \left(\phi\frac{S_r}{K_w} + \phi\frac{\partial S_r}{\partial p_w}\right)^{k+1}\dot{p}_w^{k+1} + S_r^{k+1}\mathscr{L} : \mathbf{1} + \int_{\mathscr{H}} \underline{\varrho}\left(\underline{\mathscr{Q}} - \underline{\mathscr{Q}}'\right)^{k+1} \, \mathrm{d}\mathscr{V}'. \tag{116}$$

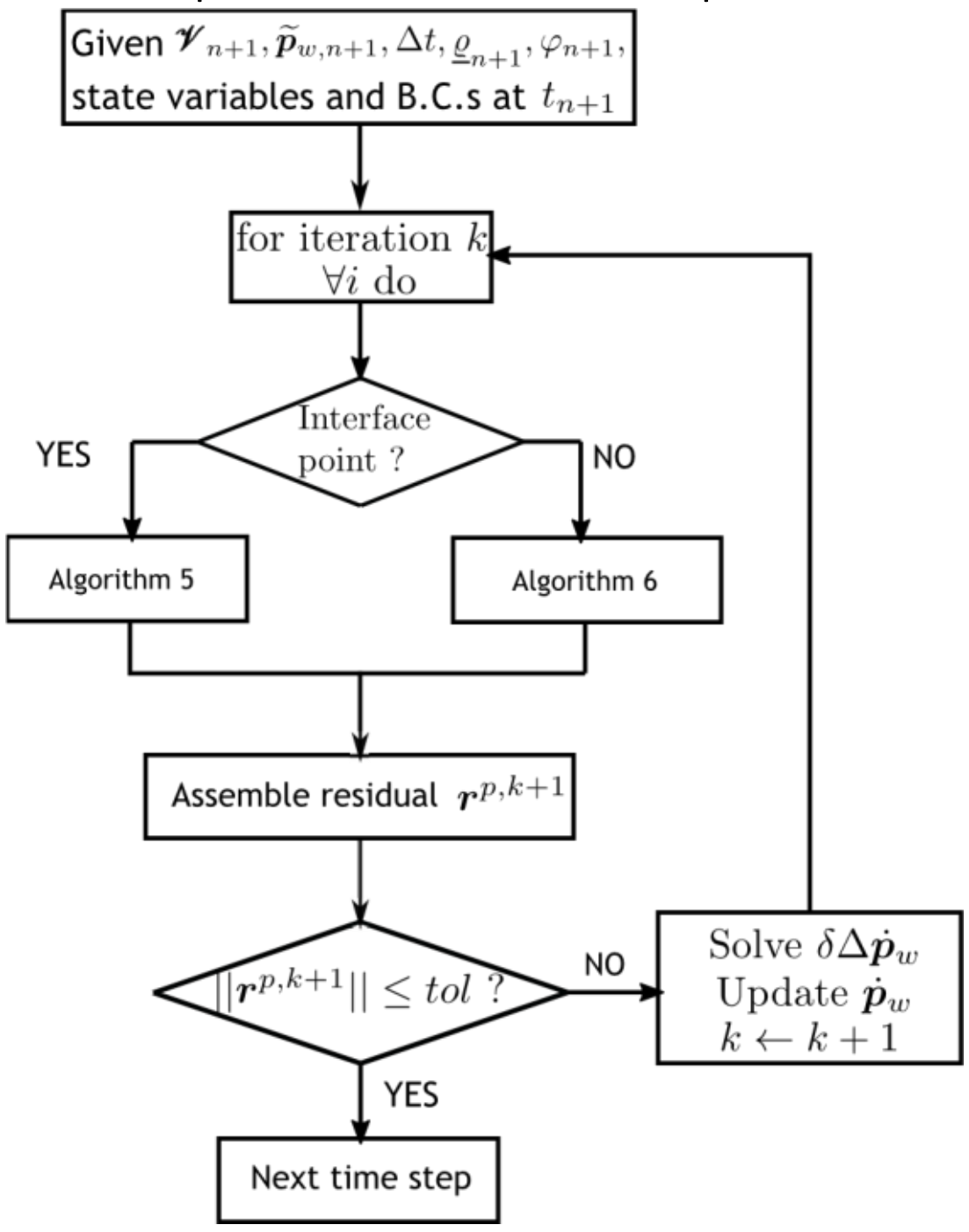

Figure 4: Flowchart of Newton's method for solving the mass balance equation in the unsaturated fluid flow stage of the fractional step algorithm for coupled updated Lagrangian-Eulerian periporomechanics.

Similarly, the residual of the mass balance equation for a phreatic material points reads

$$r^{p,k+1} = \left(\phi\frac{S_r}{K_w} + \phi\frac{\partial S_r}{\partial p_w}\right)^{k+1} \dot{p}_w^{k+1} + +S_r^{k+1}\left(\mathscr{L}^{(1)} + \mathscr{L}^{(2)}\right) : \mathbf{1}$$
$$+ \int_{\mathscr{H}} \underline{\varrho}\left\{\left(\underline{\mathscr{Q}}^{(1)} + \underline{\mathscr{Q}}^{(2)}\right) - \left(\underline{\mathscr{Q}}'^{(1)} + \underline{\mathscr{Q}}'^{(2)}\right)\right\}^{k+1} \mathrm{d}\mathscr{V}'. \tag{117}$$

Algorithms 5 and 6 show the algorithms for computing the fluid flow and volume rate states for the residual vector of the mass balance equations for the bulk and phreatic material points. Through an implicit Newton's method, $\Delta\, \dot{} p_{w,n+1}$ can be solved as follows.

$$\mathscr{R}^{p,k+1} = \mathscr{R}^{p,k} + \mathscr{A}^{u,k}\delta\Delta\dot{\boldsymbol{p}}_w^{k+1} \approx \mathbf{0}, \tag{118}$$

where $\mathbf{R}^p$ is the global residual vector of the mass balance and $\mathbf{A}^p$ is the corresponding global tangent operator

$$\mathscr{A}^{p,k} = \frac{\partial \mathscr{R}^{p,k}}{\partial \Delta\dot{\boldsymbol{p}}_w^k}. \tag{119}$$

**Algorithm 4** Fractional step algorithm for fluid flow in fully coupled updated Lagrangian periporomechanics

```
 1: procedure GIVEN 𝒱̄_{n+1}, p̃_{w,n+1}, p̃̇_{w,n+1}, p̃_{f,n+1}, p̃̇_{f,n+1}, t_n AND Δt, SOLVE ṗ_{w,n+1}
 2:     Compute fluid flow and volume change via Algorithm 5
 3:     Compute the balance of mass residual r^{p,0}_{n+1}
 4:     Set k = 0
 5:     while ||r^{p,k}|| > tol_2 do
 6:         Construct tangent 𝒜^{pf,k} = [{∂ℛ^{pf,k}/∂Δp^k_w}; {∂ℛ^{pf,k}/∂Δp^k_f}]
 7:         Solve 𝒜^{pf,k}[δΔp_w; δΔp_f]^{k+1} = −{ℛ^{p,k}, ℛ^{f,k}} for [δΔp_w; δΔp_f]^{k+1}
 8:         Update p^{k+1}_w and p^{k+1}_f
 9:         Check the residual ℛ^{pf,k+1}
10:         Set k ← k + 1
11:     end while
12:     Update u_n ← u_{n+1}, x_n ← x_{n+1},
13: end procedure
```

Solving (118) gives,

$$\delta\Delta\dot{\boldsymbol{p}}_w^{k+1} = -(\boldsymbol{\mathscr{A}}^k)^{-1}\mathscr{R}^{p,k}. \tag{120}$$

Finally, we have

$$\Delta\dot{\boldsymbol{p}}_w^{k+1} = \Delta\dot{\boldsymbol{p}}_w^k + \delta\Delta\dot{\boldsymbol{p}}_w^{k+1}. \tag{121}$$

*3.3. Spatial discretization of the updated Lagrangian-Eulerian periporomechanics*

The mixed Lagrangian-Eulerian meshfree method is adopted to spatially discretize the coupled governing equations [23]. Let P denote the number of total material points in the problem domain and $N_i$ be the number of material points in the horizon of material point *i*. The spatially discretization of the motion equation and mass balance equation can be written as

$$\mathscr{A}_{i=1}^{\mathscr{P}}(\boldsymbol{\mathscr{M}}_i\ddot{\mathbf{u}}_i - \overline{\mathscr{T}}_i + \mathscr{T}_{w,i} + \mathscr{F}_i) = 0, \tag{122}$$
$$\mathscr{A}_{i=1}^{\mathscr{P}}(\mathscr{S}_i + \mathscr{Q}_i + \overline{\mathscr{V}}_{s,i} + \mathscr{F}_i) = 0. \tag{123}$$

**Algorithm 5** Construct $\boldsymbol{\mathcal{Q}}_{n+1}$ for interface points

```
for all interface points i do
    if p_w < 0 then
        Compute degree saturation S^k_{r,i}
        Compute relative permeability k_i^{r,k}
    end if
    if p_w ≥ 0 then
        Set S^k_{r,i} = 1
        Set k_i^{r,k} = 1
    end if
    for all neighbors j do
        if α_ij = 0 then
            Compute pressure gradient grad(Φ̃)_i^{k,(1)}
        end if
        if α_ij = 1 then
            Compute pressure gradient grad(Φ̃)_i^{k,(2)}
        end if
    end for
    Compute the flux vector q_i^{k,(1)} and q_i^{k,(2)}
end for
for all interface points i do
    for all neighbors j do
        if α_ij = 0 then
            Compute Q_ij^{k,(1)}
        end if
        if α_ij = 1 then
            Compute Q_ij^{k,(2)}
        end if
    end for
end for
```

Where

$$\boldsymbol{\mathcal{M}}_i = [\rho_s(1-\phi_i) + \rho_w S_{r,i}\phi_i]\,\mathcal{V}_i \mathbf{1}, \tag{124}$$

$$\overline{\boldsymbol{\mathcal{T}}}_i = \sum_{j=1}^{\mathcal{N}_i} \left[ \left( \underline{\overline{\boldsymbol{\mathcal{T}}}}_{ij}^{(1)} + \underline{\overline{\boldsymbol{\mathcal{T}}}}_{ij}^{(2)} \right) - \left( \underline{\overline{\boldsymbol{\mathcal{T}}}}_{ji}^{(1)} + \underline{\overline{\boldsymbol{\mathcal{T}}}}_{ji}^{(2)} \right) \right] \mathcal{V}_j \mathcal{V}_i, \tag{125}$$

$$\boldsymbol{\mathcal{T}}_{w,i} = \sum_{j=1}^{\mathcal{N}_i} \left[ S_{r,i} \left( \underline{\boldsymbol{\mathcal{T}}}_{w,ij}^{(1)} + \underline{\boldsymbol{\mathcal{T}}}_{w,ij}^{(2)} \right) - S_{r,j} \left( \underline{\boldsymbol{\mathcal{T}}}_{w,ji}^{(1)} + \underline{\boldsymbol{\mathcal{T}}}_{w,ji}^{(2)} \right) \right] \mathcal{V}_j \mathcal{V}_i, \tag{126}$$

$$\mathcal{S}_i = \phi_i \left( \frac{S_r}{K_w} - \frac{\partial S_{r,i}}{\partial s_i} \right) \dot{p}_{w,i} \mathcal{V}_i, \tag{127}$$

*A* is an global assembly operator [23, 54], and $V_i$ and $V_j$ are the volumes of material points *i* and *j*, respectively, in the current configuration. It is noted that (125), (126), (128) (129) and are written in the forms for the phreatic material points in which only one term remains for a bulk material points.

In (125) and (126), the effective force state and the water force state are written as

$$\underline{\overline{\mathscr{T}}}^{(1)}_{ij} = \left\{\underline{\omega}_{ij}\underline{\boldsymbol{\zeta}}_{ij}\overline{\boldsymbol{\sigma}}_i\right\}^{(1)} \boldsymbol{\mathscr{K}}_i^{-1} + \left\{\underline{\omega}_{ij}\frac{GC}{\omega_0}\boldsymbol{\mathscr{R}}_i^s\right\}^{(1)}, \tag{130}$$

$$\underline{\overline{\mathscr{T}}}^{(2)}_{ij} = \left\{\underline{\omega}_{ij}\underline{\boldsymbol{\zeta}}_{ij}\overline{\boldsymbol{\sigma}}_i\right\}^{(2)} \boldsymbol{\mathscr{K}}_i^{-1} + \left\{\underline{\omega}_{ij}\frac{GC}{\omega_0}\boldsymbol{\mathscr{R}}_i^s\right\}^{(2)}, \tag{131}$$

$$\underline{\mathscr{T}}^{(1)}_{w,ij} = \left\{\underline{\omega}_{ij}\underline{\boldsymbol{\zeta}}_{ij}\mathbf{1}p_{w,i}\right\}^{(1)} \boldsymbol{\mathscr{K}}_i^{-1}, \tag{132}$$

$$\underline{\mathscr{T}}^{(2)}_{w,ij} = \left\{\underline{\omega}_{ij}\underline{\boldsymbol{\zeta}}_{ij}\mathbf{1}p_{w,i}\right\}^{(2)} \boldsymbol{\mathscr{K}}_i^{-1}. \tag{133}$$

The fluid flow states in (128) read

$$\underline{\mathscr{Q}}^{(1)}_{ij} = \left\{\underline{\omega}_{ij}\rho_w\underline{\boldsymbol{\zeta}}_{ij}\boldsymbol{q}_i\right\}^{(1)} \boldsymbol{\mathscr{K}}_i^{-1} + \left\{\underline{\omega}_{ij}\frac{GK_p}{\omega_0}\mathscr{R}_i^w\right\}^{(1)}, \tag{134}$$

$$\underline{\mathscr{Q}}^{(2)}_{ij} = \left\{\underline{\omega}_{ij}\rho_w\underline{\boldsymbol{\zeta}}_{ij}\boldsymbol{q}_i\right\}^{(2)} \boldsymbol{\mathscr{K}}_i^{-1} + \left\{\underline{\omega}_{ij}\frac{GK_p}{\omega_0}\mathscr{R}_i^w\right\}^{(2)}. \tag{135}$$

The velocity gradients in (129) read

$$\boldsymbol{\mathscr{L}}_{ij,n+1} = \left[\sum_{k=1}^{\mathscr{N}_i}\left(\underline{\omega}_{ij}\dot{\underline{\boldsymbol{\mathscr{Y}}}}_{ik}\otimes\underline{\boldsymbol{\zeta}}_{ik}\right)\mathscr{V}_k\right]^{(1)} (\boldsymbol{\mathscr{K}}_{ij})^{-1}, \tag{136}$$

$$\boldsymbol{\mathscr{L}}_{ij,n+1} = \left[\sum_{k=1}^{\mathscr{N}_i}\left(\underline{\omega}_{ij}\dot{\underline{\boldsymbol{\mathscr{Y}}}}_{ik}\otimes\underline{\boldsymbol{\zeta}}_{ik}\right)\mathscr{V}_k\right]^{(2)} (\boldsymbol{\mathscr{K}}_{ij})^{-1}. \tag{137}$$

In the next section, we present numerical examples to validate the computational implementation and demonstrate the efficacy of the proposed computational periporomechanics for modeling dynamic failure and fracturing in unsaturated porous media.

**Algorithm 6** Construct $\boldsymbol{\mathcal{Q}}_{n+1}$ for non-interface points

**for** all points $i$ **do**
    **if** $p_{w,(i)} < 0$ **then**
        Compute degree saturation $S^k_{r,i}$
        Compute relative permeability $k^{r,k}_i$
    **end if**
    **if** $p_{w,i} \geq 0$ **then**
        Set $S^k_{r,i} = 1$
        Set $k^{r,k}_i = 1$
    **end if**
    Compute pressure gradient $\mathrm{grad}(\widetilde{\Phi})^k_i$
    Compute the flux vector $\boldsymbol{q}^k_i$
**end for**
**for** all points $i$ **do**
    **for** all neighbors $j$ **do**
        Compute fluid flow states $\underline{\mathscr{Q}}^k_{ij}$
    **end for**
**end for**

## 4. Numerical examples

### *4.1. Example 1: Validation of the coupled stabilization scheme*

In this example, we simulate the dynamic consolidation of a three-dimensional unsaturated soil specimen to validate the proposed coupled stabilization scheme. Figure 5 presents the geometry of the soil column, loading protocol, and boundary conditions. The problem domain is discretized into 9,000 uniform mixed material points. The distance of two neighboring material point centers is $\Delta x$ = 0.1 m. For the fluid phase, the bottom boundary is prescribed constant fluid pressure, and all other boundaries are impervious. The solid phase is modeled using an isotropic elastic correspondence model [22, 23]. The intrinsic permeability is assumed isotropic and uniform. The material properties chosen are $\rho_s = 2.1\times10^3$ kg/m$^3$, $\rho_w = 1\times10^3$ kg/m$^3$, $\mu_w = 1\times10^{-3}$ Pa·s, initial porosity $\varphi_0 = 0.33$, bulk modulus $K = 3.3\times10^4$ kPa, shear modulus $\mu_s = 1.62 \times 10^4$ kPa, water bulk modulus $K_w = 2 \times10^5$ kPa, $k_w = 1 \times 10^{-14}$ m$^2$, $n = 1.25$, $s_a = 10$ kPa. The initial uniform effective stress $\sigma$ = -12.33 kPa. The initial fluid pressure $p_w$ = -15 kPa (i.e., $S_r = 0.82$). The loading rate is $\dot{u}_y = 0.01$ m/s. The total loading time $t$ = 5 s with the time increment $\Delta t$ = 0.005 s.

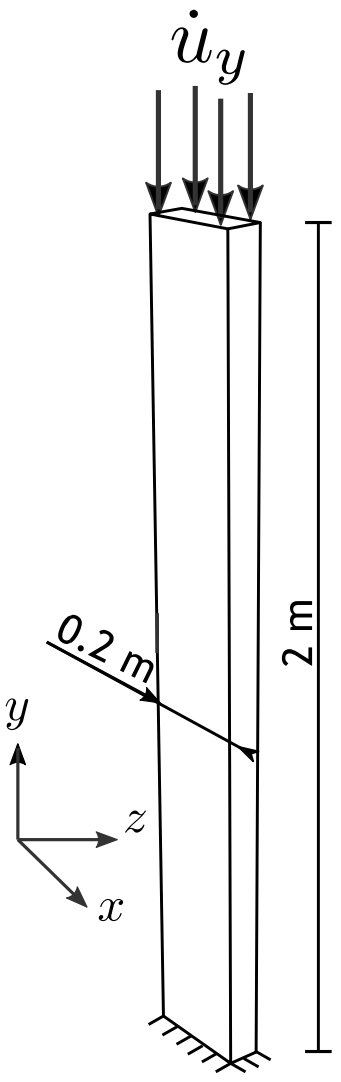

Figure 5: Problem setup for example 1.

Figures 6 and 7 compare the contours of displacement and pressure, respectively, from the simulations using different values of $G$. The results of simulations with the standard correspondence material model ($G$ = 0) show noticeable oscillations. However, the oscillations have disappeared in the results with the stabilized correspondence material model ($G$ = 1.0).

Figure 8 plots the results of the vertical velocity within the specimen with three values of $G$ at $u_y$ = 0.05 m on the top boundary. Figure 9 plots the fluid pressure within the specimen with three values of $G$ at the same load step. The results in Figures 8 and 9 show the effect of stabilization could depend on the value of $G$. This preliminary study implies that the simulation with $G$ = 1.0 removes the oscillations of the vertical velocity and fluid pressure in this example.

### *4.2. Example 2: Variably saturated flow with an evolving phreatic interface in porous media*

In this example, we simulate variably saturated fluid flow with the evolution of the phreatic interface in a rigid porous material. Figure 10 presents the problem setup and the contour of the initial fluid pressure with the phreatic line. The initial fluid pressure is prescribed through a linear function of the distance above and below the phreatic line as shown in Figure 10. The solid skeleton is assumed rigid. All fluid boundaries are impervious. The problem domain is discretized into 15,000 uniform mixed material points with $\Delta x$ = 0.15 m. The initial phreatic line is determined through the search algorithm 2 introduced in Section 3. The material properties chosen are $\rho_w = 1 \times10^3$ kg/m$^3$, $\mu_w = 1 \times 10^{-3}$ Pa.s, $K_w = 2 \times 10^5$ kPa, $k_w = 1 \times 10^{-15}$ m$^2$, $a_1 = 0.038$, $a_2 = 3.49$, $n = 3.0$, $s_a = 50$ kPa. The horizon is $\delta = 3.05\Delta x$. The total simulation time $t$ = 8 h (hour) with $\Delta t$ = 1 s.

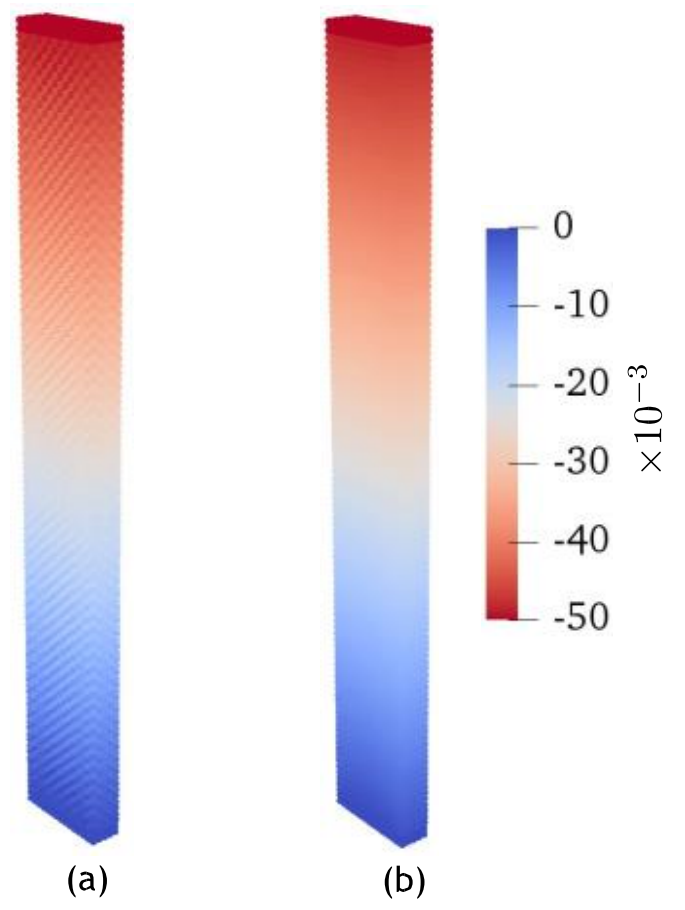

Figure 6: Contours of the vertical displacement with (a) $G$ = 0.0 and (b) $G$ = 1.0 at $u_y$ = 0.05 m.

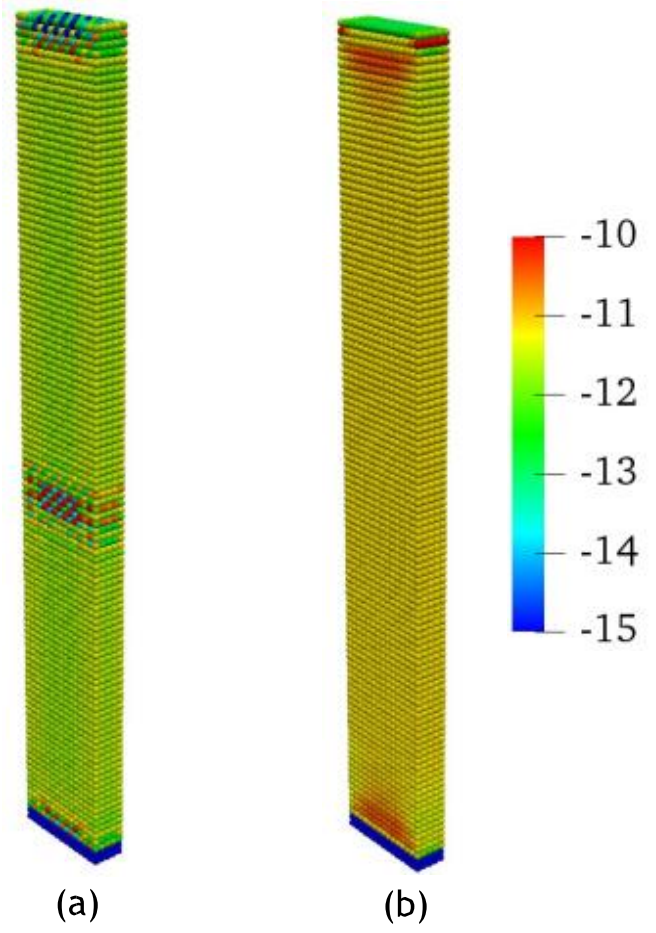

Figure 7: Contours of water pressure (kPa) with (a) $G$ = 0.0 and (b) $G$ = 1.0 at $u_y$ = 0.05 m.

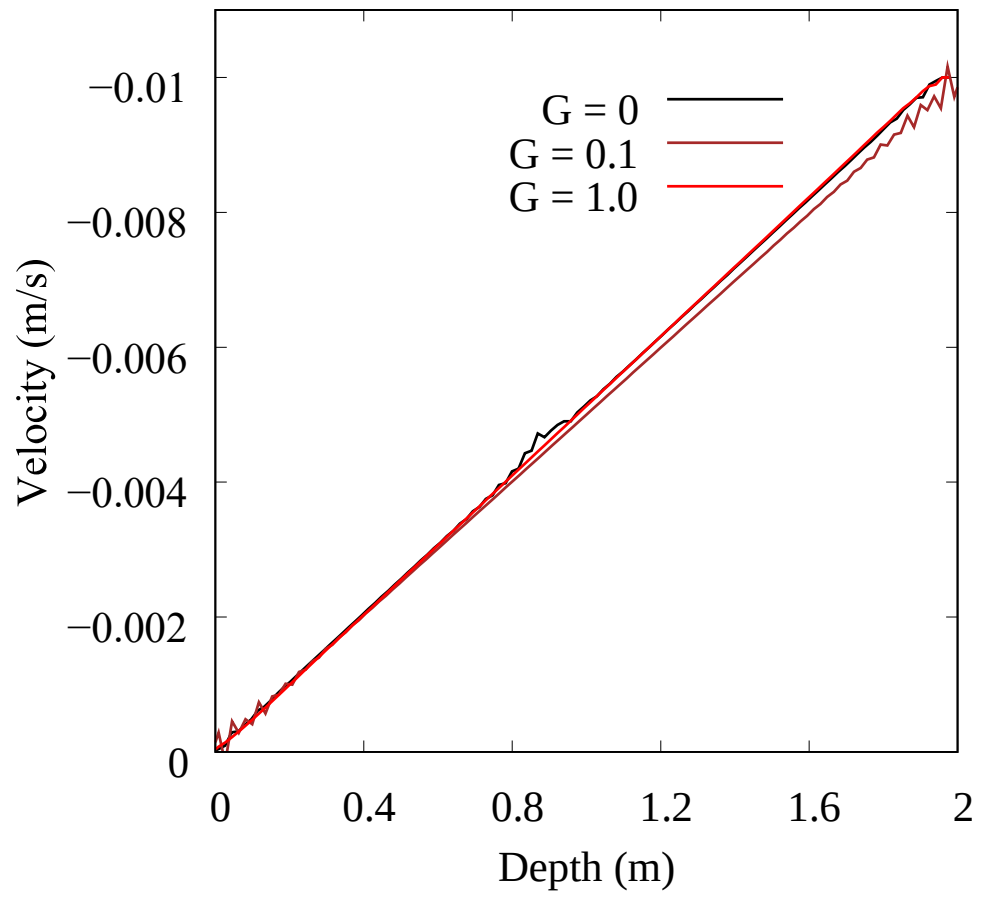

Figure 8: Variation of the vertical velocity along the specimen at $u_y$ = 0.05 m from the simulations with three values of $G$.

Figures 11 and 12 plot the results of the simulation at (a) $t$ = 1 h, (b) $t$ = 4 h, and (c) $t$ = 8 h. Fluid flow in the problem domain is driven by initial condition of the pressure in the domain. The initial spatial variation of fluid pressure drives fluid flow toward the unsaturated zone causing the phreatic interface to move toward the right bottom of the problem domain. The initial straight phreatic line becomes curved because of the saturated and unsaturated fluid flow within the problem domain. This example demonstrates that the implemented variably saturated periporomechanics with a dedicated search algorithm for phreatic material points can adequately model the evolution of the phreatic interface under partially saturated fluid flow condition in porous media.

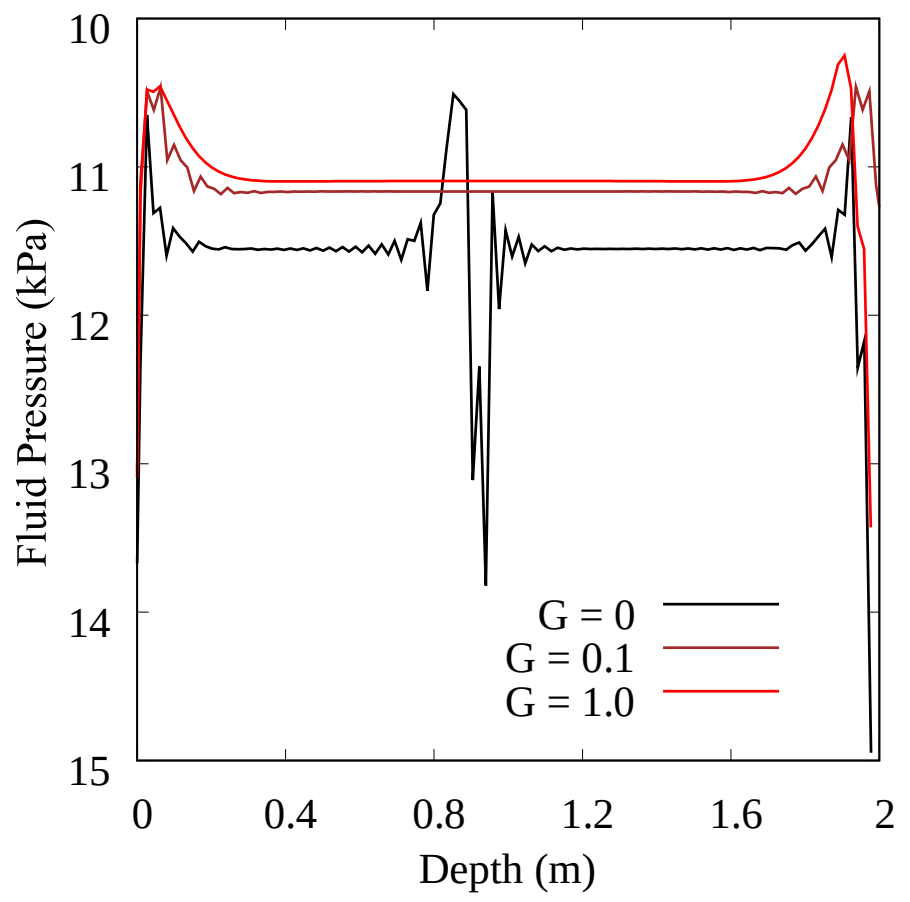

Figure 9: Variation of water pressure along the specimen at $u_y$ = 0.05 m from the simulations with three values of *G*.

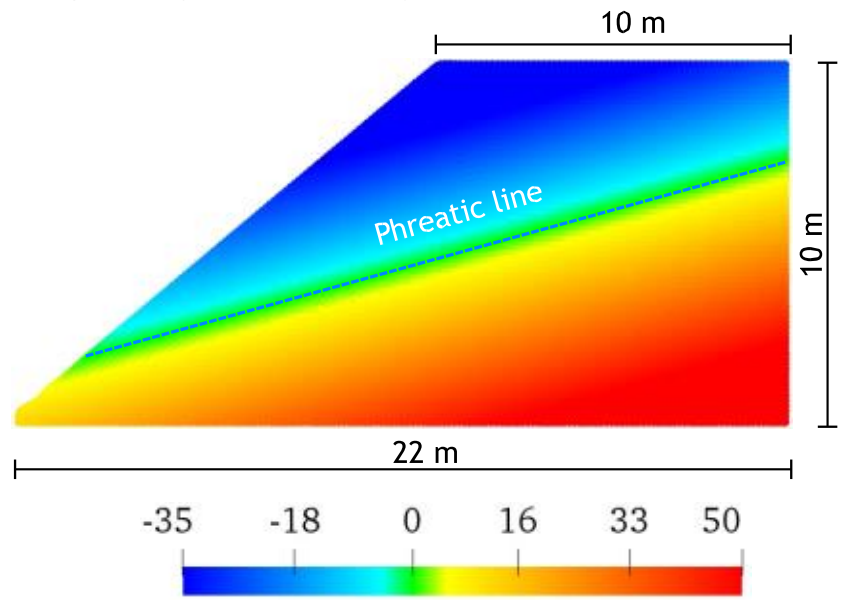

Figure 10: Problem setup for example 2 and the contour of initial fluid pressure (kPa) with a phreatic line.

### *4.3. Example 3: Mode I cracking*

In this example, we simulate mode I crack propagation and unsaturated fluid flow in a porous body in two dimensions. We compare the numerical results from the updatedLagrangian formulation in this study to the results from the total-Lagrangian periporomechanics in [24]. Figure 13 presents the problem geometry, boundary conditions, and loading protocol. The pre-existing crack is defined by eliminating interaction between material points across the crack plane [24]. All fluid phase boundaries of the specimen are assumed impermeable. The problem domain is discretized into 20000 uniform mixed material points with $\Delta x$ = 0.0025 m.

The solid phase is modeled using an isotropic elastic correspondence constitutive model [22, 23]. The material properties are $\rho_s$ = 2×10$^3$ kg/m$^3$, $\rho_w$ = 1×10$^3$ kg/m$^3$, $\mu_w$ = 1×10$^{-3}$ Pa·s, initial porosity $\varphi_0$ = 0.25, $G_0$ = 225 J/m$^2$, $K$ = 1.346 × 10$^7$ kPa, $\mu_s$ = 1.095 × 10$^7$ kPa, intrinsic permeability $k_w$ = 1×10$^{-16}$ m$^2$, $a_1$ = 0.038, $a_2$ = 3.49, $n$ = 1.78, $s_a$ = 1.2×10$^4$ kPa. The horizon $\delta$ = 3.05 $\Delta x$. The stabilization parameter $G$ = 1.0. The initial uniform effective stress $\bar{\sigma}$ = -13 kPa. The initial water pressure in the specimen is -15 kPa ($S_r$ = 0.87). The loading rate $\dot{u}_y$ = 2.0 ×10$^{-6}$ m/s. The total loading time $t$ = 1000 s with $\Delta t$ = 0.5 s.

The results are presented in Figures 14, 15, and 16. Figure 14 compares the reaction force versus the applied displacement on the top boundary from the simulations using the updated Lagrangian formulation to the results from the total-Lagrangian formulation. The results in Figure 14 show good agreement between the two curves before the peak load. The slight difference between the two curves after the peak load could be due to the fact that the neighboring points of a mixed material point are updated on each time step in the updated Lagrangian formulation.

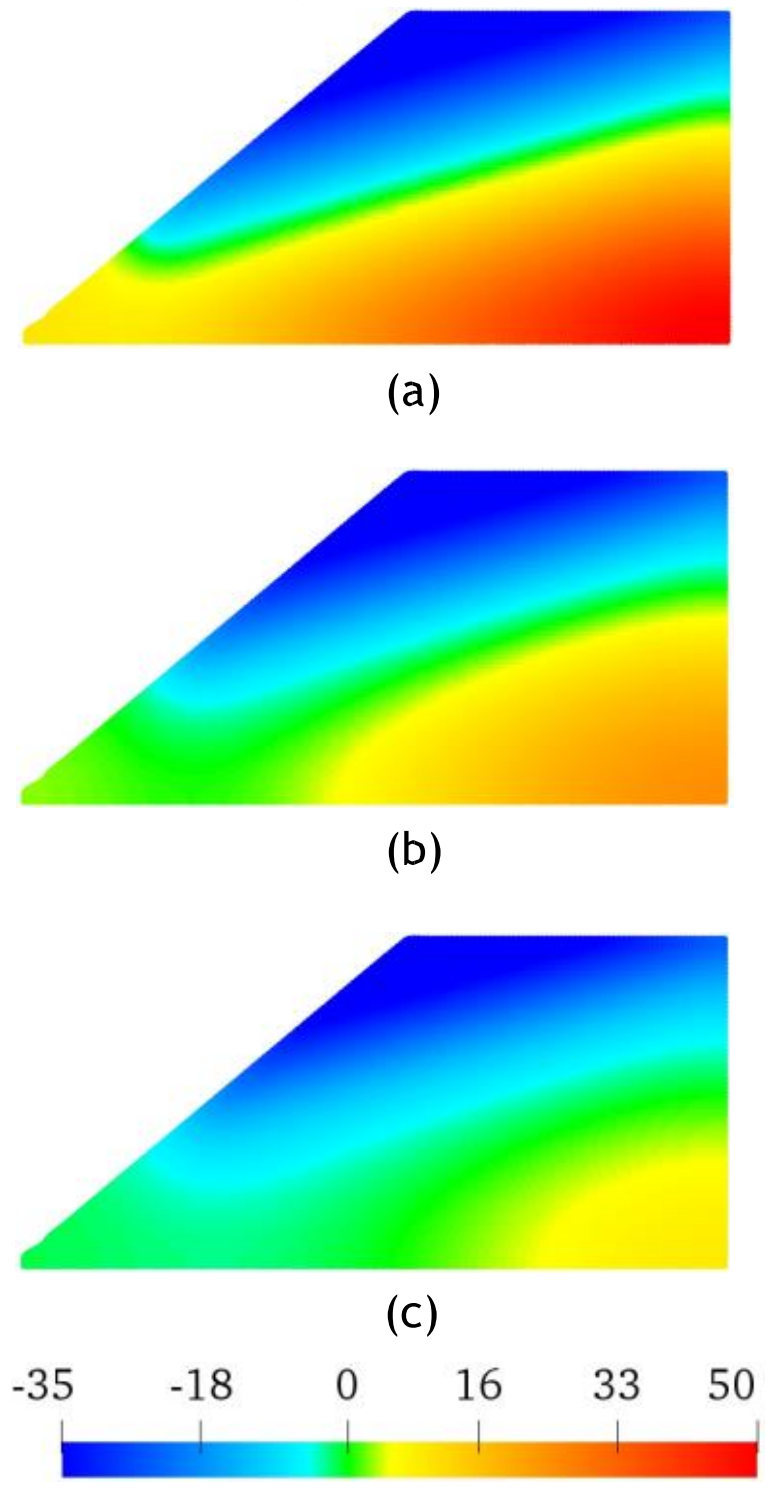

Figure 11: Contours of the water pressure (kPa) in the problem domain at (a) $t$ = 1 h, (b) $t$ = 4 h and (c) $t$ = 8 h.

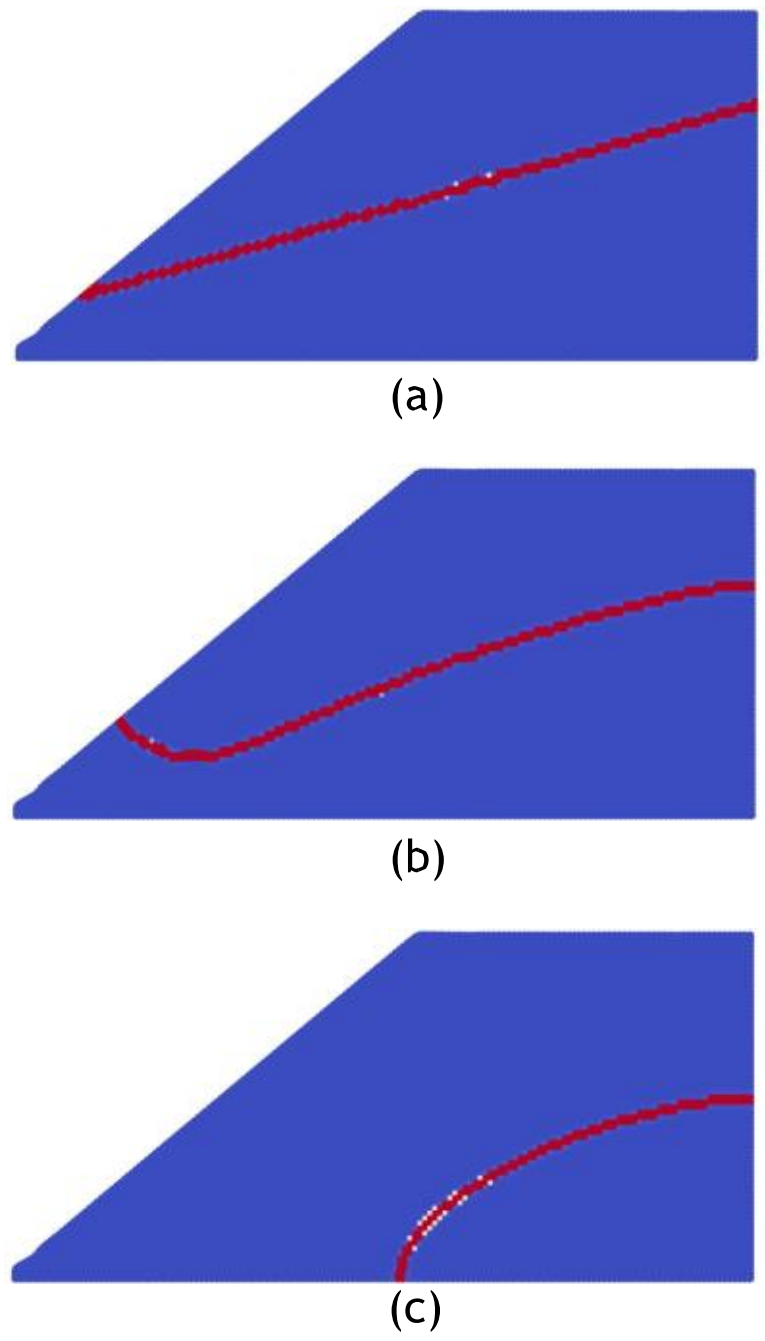

Figure 12: Evolution of the phreatic interface (red line) in the problem domain at (a) $t$ = 1 h, (b) $t$ = 4 h and (c) $t$ = 8 h.

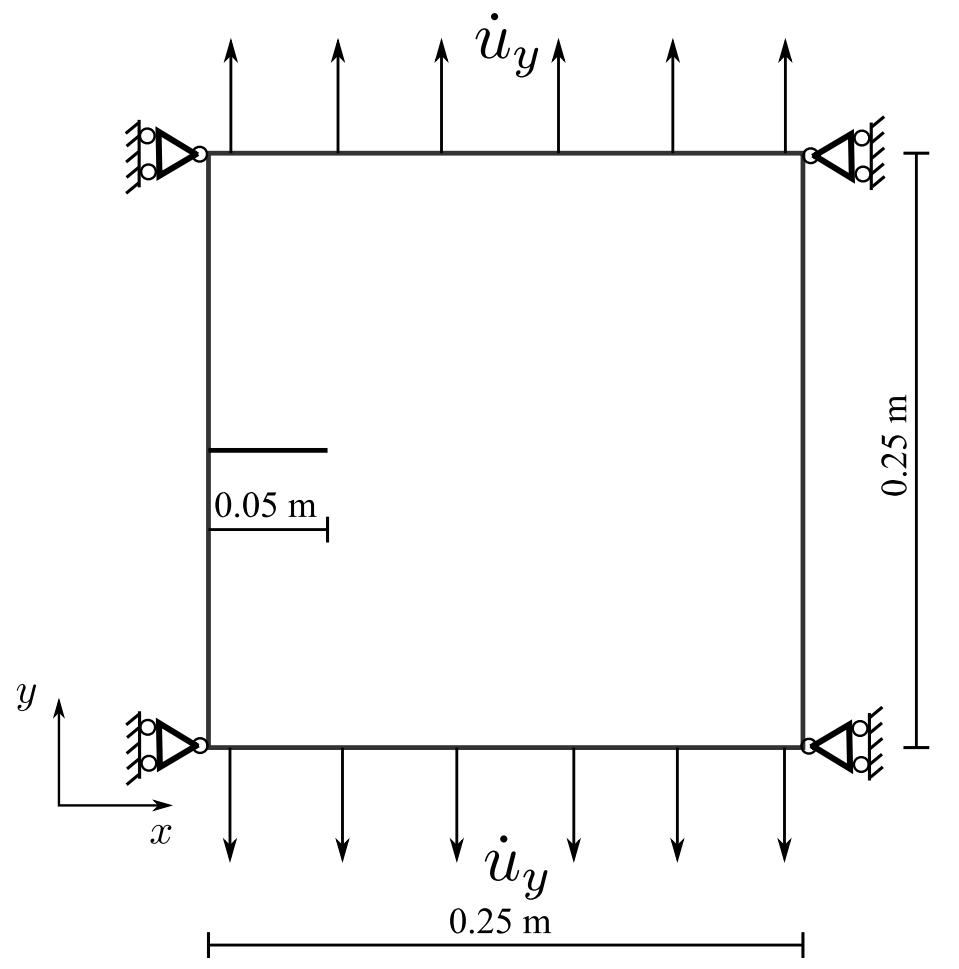

Figure 13: Problem setup for example 3.

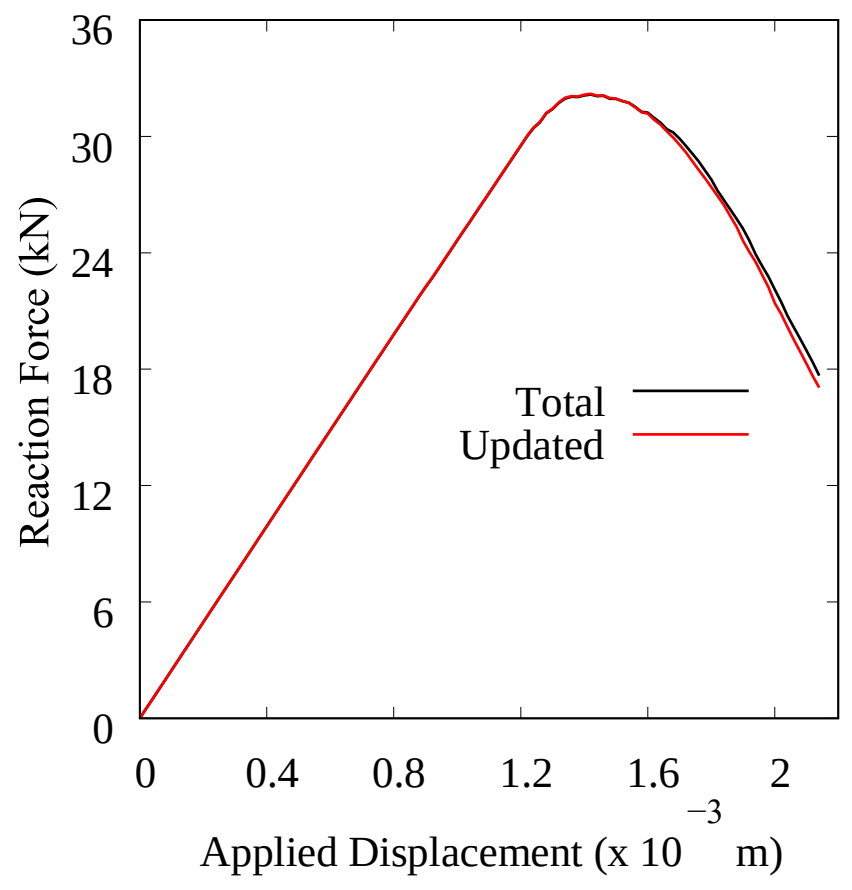

Figure 14: Comparison of the loading curves from the total and updated Lagrangian formulations.

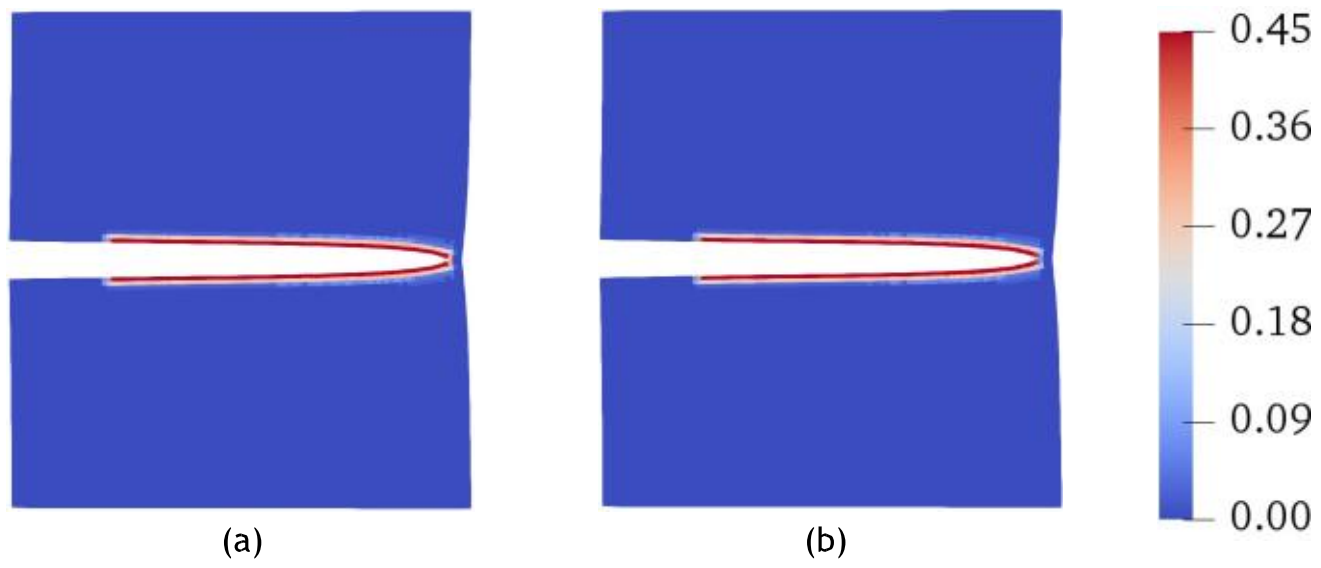

Figure 15: Contours of damage variable ($\phi$) from the simulations using (a) total and (b) updated Lagrangian formulations at $u_y$ = 0.002 m (×50).

Figures 15 and 16 plot the contours of damage variable and water pressure at $u_y$ = 0.002 m respectively from the simulations using total and updated Lagrangian formulations. The results in both figures show that the crack opening leads to a similar development of negative pressures (larger matric suction) in the vicinity of the crack. The negative water pressures are due to the fluid flow into the crack space in the mode I cracking process. The results from the total and updated formulation are in good agreement, while there is a slight difference in both the value of water pressure and the length of crack propagation. We note that the mode I crack propagation in this example is brittle in nature. Thus, it is postulated that large deformation is not required for the crack to propagate. Indeed, this postulation is supported

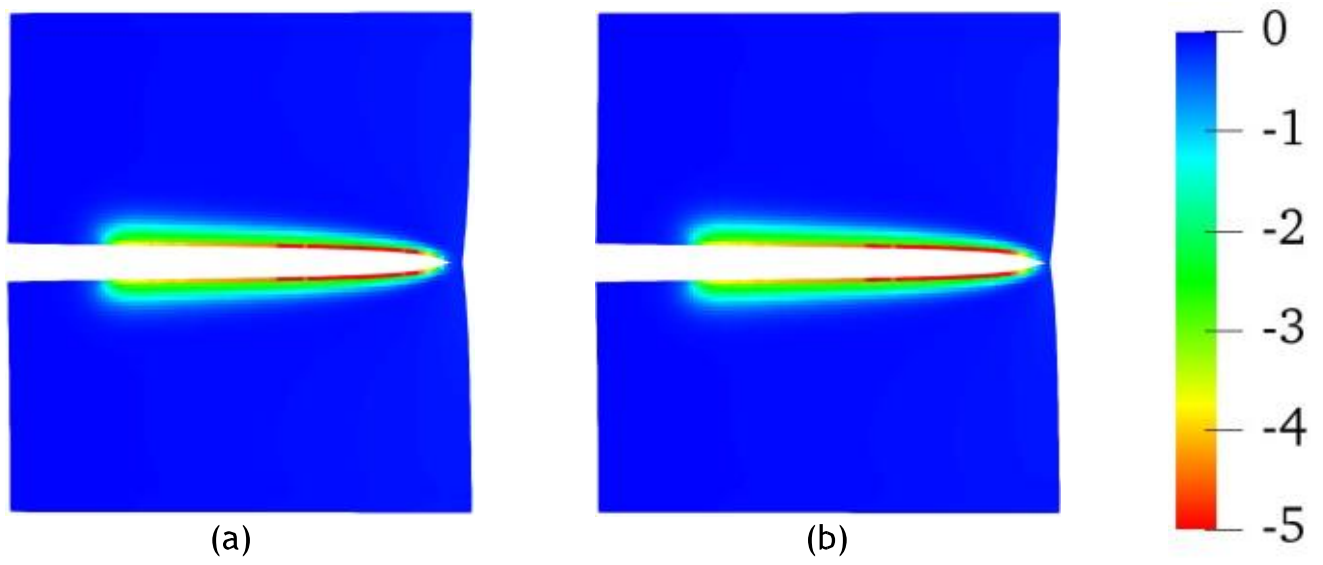

Figure 16: Contours of water pressure (MPa) from the simulations using (a) total and (b) updated Lagrangian formulations at $u_y$ = 0.002 m (×50).

by the consistency between the results from the updated and total Lagrangian formulations.

### *4.4. Wetting collapse of an unsupported vertical cut*

In this example, we simulate the wetting collapse of an unsupported vertical cut in an unsaturated soil. Figure 17 depicts the problem domain and boundary conditions. The problem domain is discretized into 13,000 mixed material points with Δ*x* = 0.1 m.

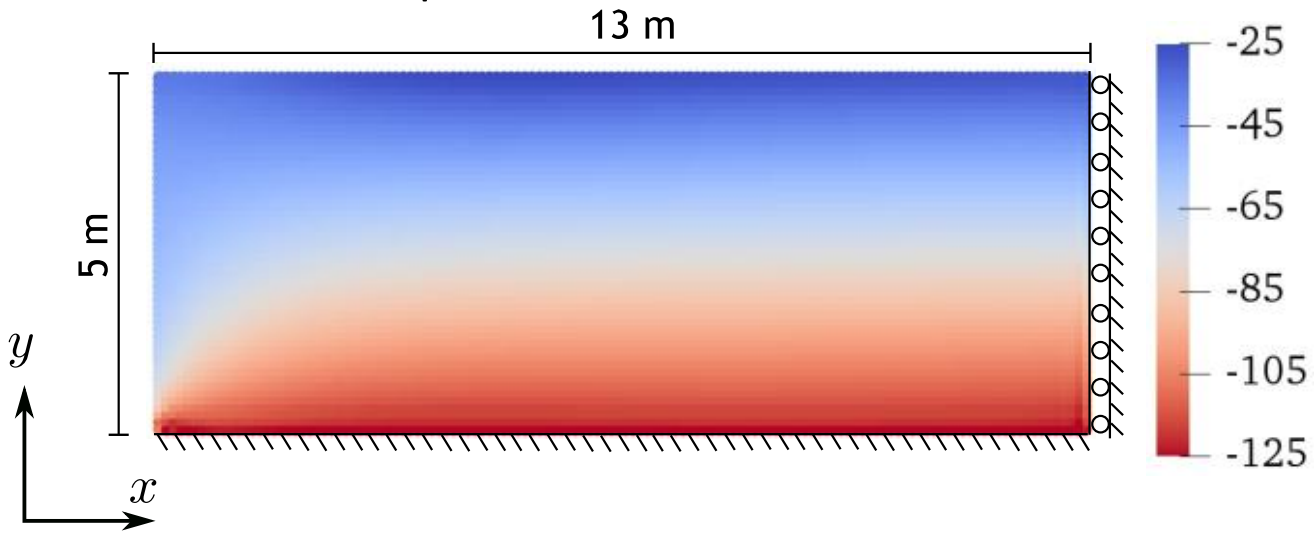

Figure 17: Problem setup for example 4 and the contour of the initial effective stress in y-direction (kPa).

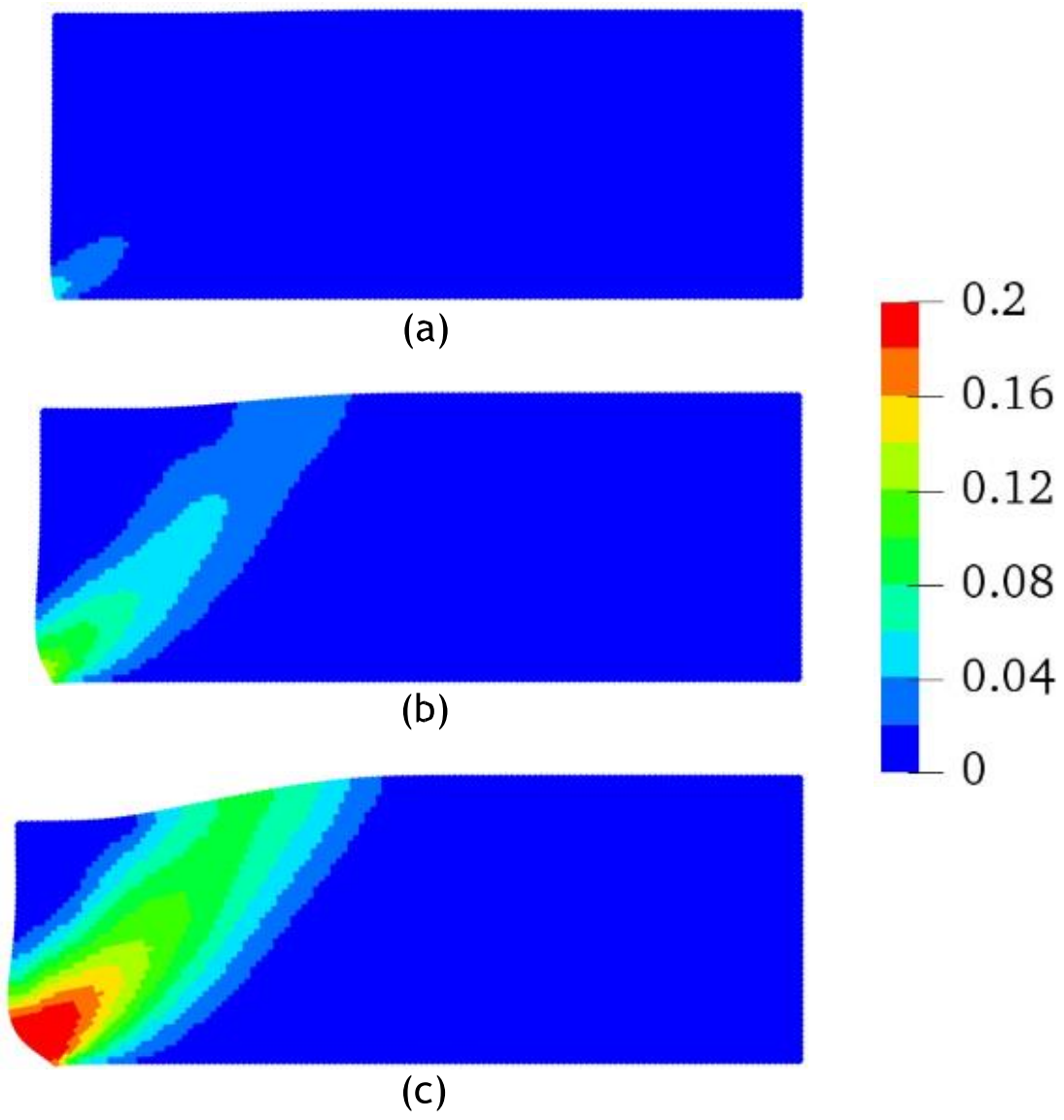

Figure 18: Contours of plastic shear strain at (a) $p_w$ = -23.3 kPa, (b) $p_w$ = -21.9 kPa and (c) $p_w$ = -20.5 kPa. (× 2)

The solid skeleton is modeled using a critical state elastoplastic constitutive model for unsaturated porous media in [23]. The material parameters are $\rho_s = 2.3 \times 10^3$ kg/m$^3$, $\rho_w = 1 \times 10^3$ kg/m$^3$, $\mu_w = 1 \times 10^{-3}$ Pa·s, $\varphi_0 = 0.25$, $K = 5.56 \times 10^4$ kPa, $\mu_s = 1.064 \times 10^4$ kPa, $K_w = 2 \times 10^5$ kPa, $k_w = 1 \times 10^{-8}$ m/s, $a_1 = 0.038$, $a_2 = 3.49$, $n = 1.78$, $s_a = 25$ kPa. For the plastic model, initial pre-consolidation pressure $\bar{p}_{c0}$ = -500 kPa, specific volume at unit pre-consolidation pressure $N = 2.2$, critical state line slope $M = 1$, swelling index $\kappa = 0.03$, compression index $\lambda = 0.13$, and fitting parameters $b_1 = 0.185$ and $b_2 = 1.42$. The horizon $\delta = 3.05\Delta x$. The stabilization parameter $G = 0.1$. The initial water pressure $p_w = -50$ kPa (i.e., matric suction = 50 kPa) throughout the problem domain. The vertical effective stress is prescribed by the equation $\sigma_y = \rho_s g h + S_r p_{w,0}$, where $h$ is the depth from the top surface. The horizontal effective stress $\sigma_x = S_r p_{w,0}$. In this example, the wetting collapse of the soil is simulated by by uniformly reducing the matric suction in the problem domain.

Figures 18 and 19 plot the contours of plastic shear and plastic volume strain, respectively, at three load steps of (a) $p_w$ = -23.3 kPa, (b) $p_w$ = -21.9 kPa, and (c) $p_w$ = -20.5 kPa. Here the plastic shear strain $\varepsilon^p_s$ is the equivalent plastic strain, i.e., $\varepsilon^p_s = {}_3$ results in 18 reveals that the failure is initiated at the left-bottom corner of the vertical cut from where the deformation band propagates upward. Figure 19 show that the maximum dilatation appears to occur at the free upper surface of the soil. As the soil deforms outward, this zone of dilatation propagates downward along with the shear band as shown in Figures 18 and 19 (b) and (c).

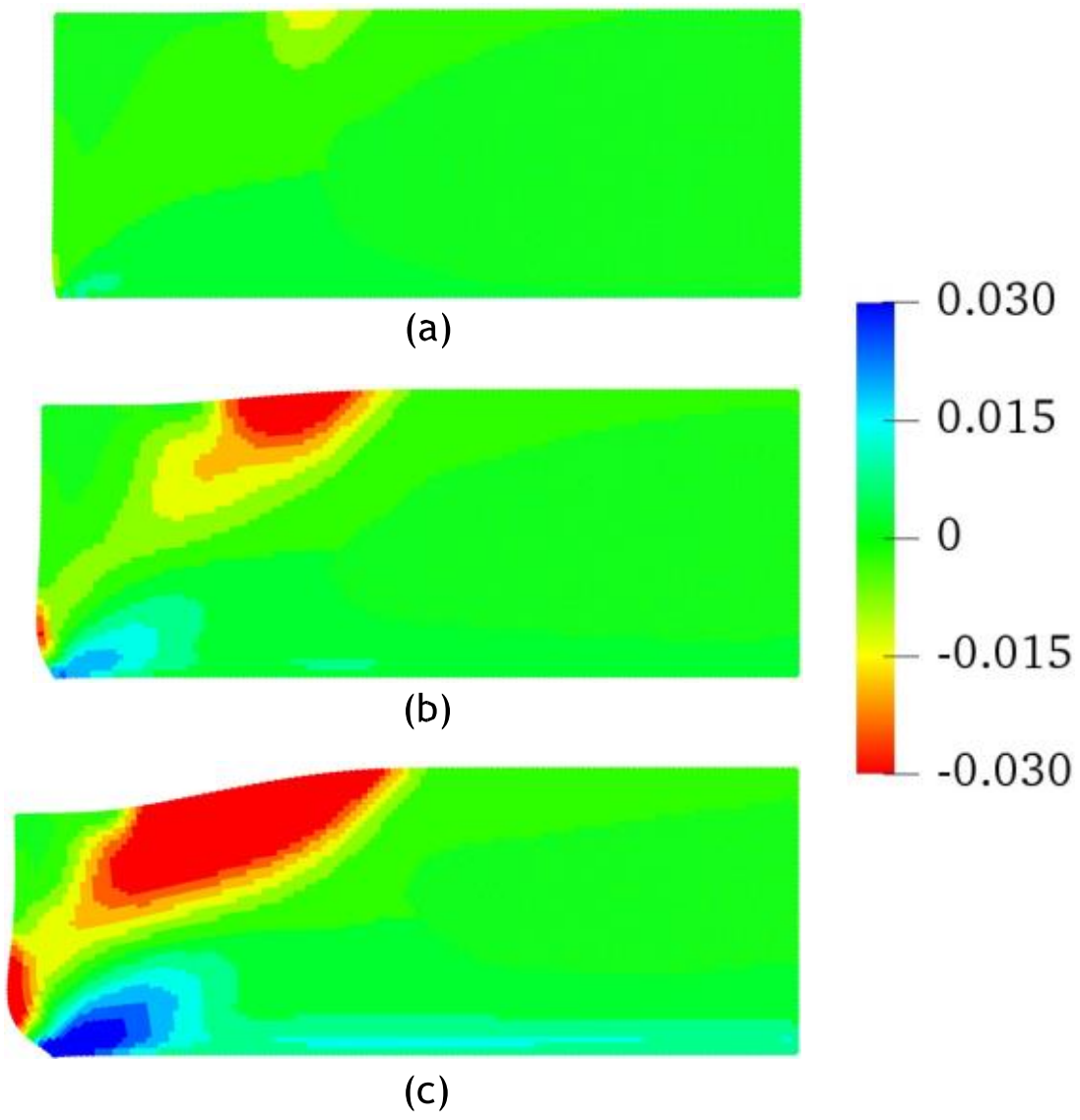

Figure 19: Contours of plastic volume strain at (a) $p_w$ = -23.3 kPa, (b) $p_w$ = -21.9 kPa, and (c) $p_w$ = -20.5 kPa. (× 2)

Next, we repeat the same simulation with a finer spatial discretization with 26,000 material points. The same horizon is utilized for both analyses. The results of the two simulations are compared in Figures 20 and 21. Figure 20 plots the contour of plastic shear strain at the same load step. Figure 21 presents the contour of plastic volume strain at the same load step. The results demonstrate that the contours of shear and volume strain from both simulations are in good agreement. The width of the deformation zone and magnitude of the plastic deformation appear insensitive to the two spatial discretizations. It may be due to the same horizon adopted for the two simulations.

### *4.5. Dynamic failure of a dam triggered by earthquake*

In this example, we demonstrate the efficacy of the proposed computational periporomechanics for modeling the dynamic failure of variably saturated porous media under earthquake loading. For this purpose, we simulate the dynamic failure of the Lower San Fernando Dam during the 1971 San Fernando earthquake ($M_w$ = 6.6). We refer to the celebrated literature [12, 55, 56] for more information about the post-earthquake investigation of this event. Figure 22 shows (a) the sketch of the geometry of the dam reported in the literature Seed et al. [56], (b) the geometry and boundary condition adopted in this study, and (c) the water table in the reservoir and the phreatic line within the dam. The base of the dam is constrained against the vertical motion, and the horizontal ends are constrained against lateral movement. The problem domain is discretized into 60,000 uniform mixed points with $\Delta x$ = 0.33 m. The horizon $\delta$ = 3.05$\Delta x$. The stabilization parameter $G$ = 0.25.

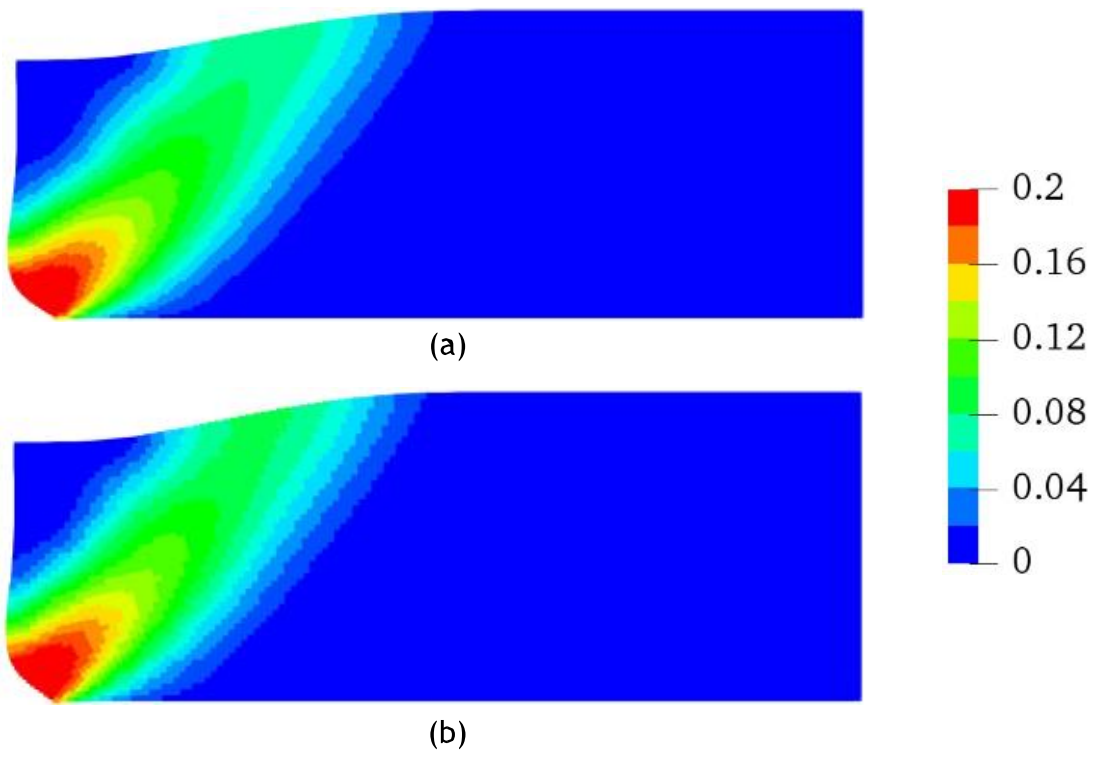

Figure 20: Contours of plastic shear strain from the simulations with (a) 26000 material points and (b) 13000 material points. (× 2)

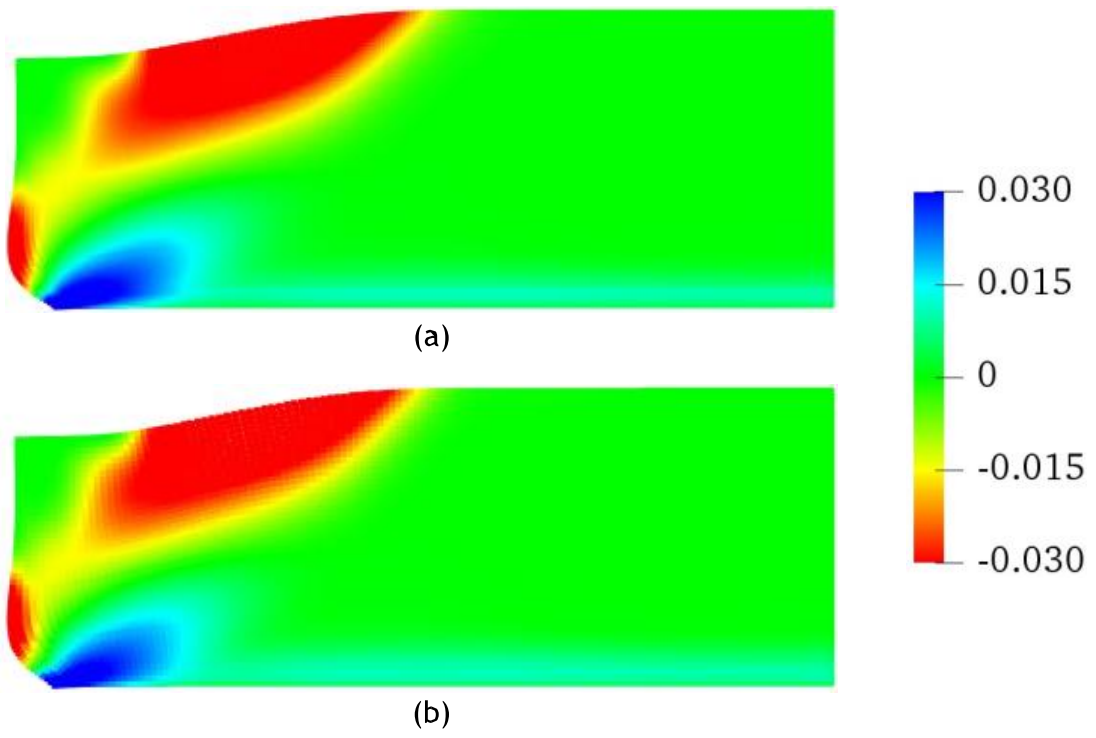

Figure 21: Contours of plastic volume strain from the simulations with (a) 26000 material points and (b) 13000 material points. (× 2)

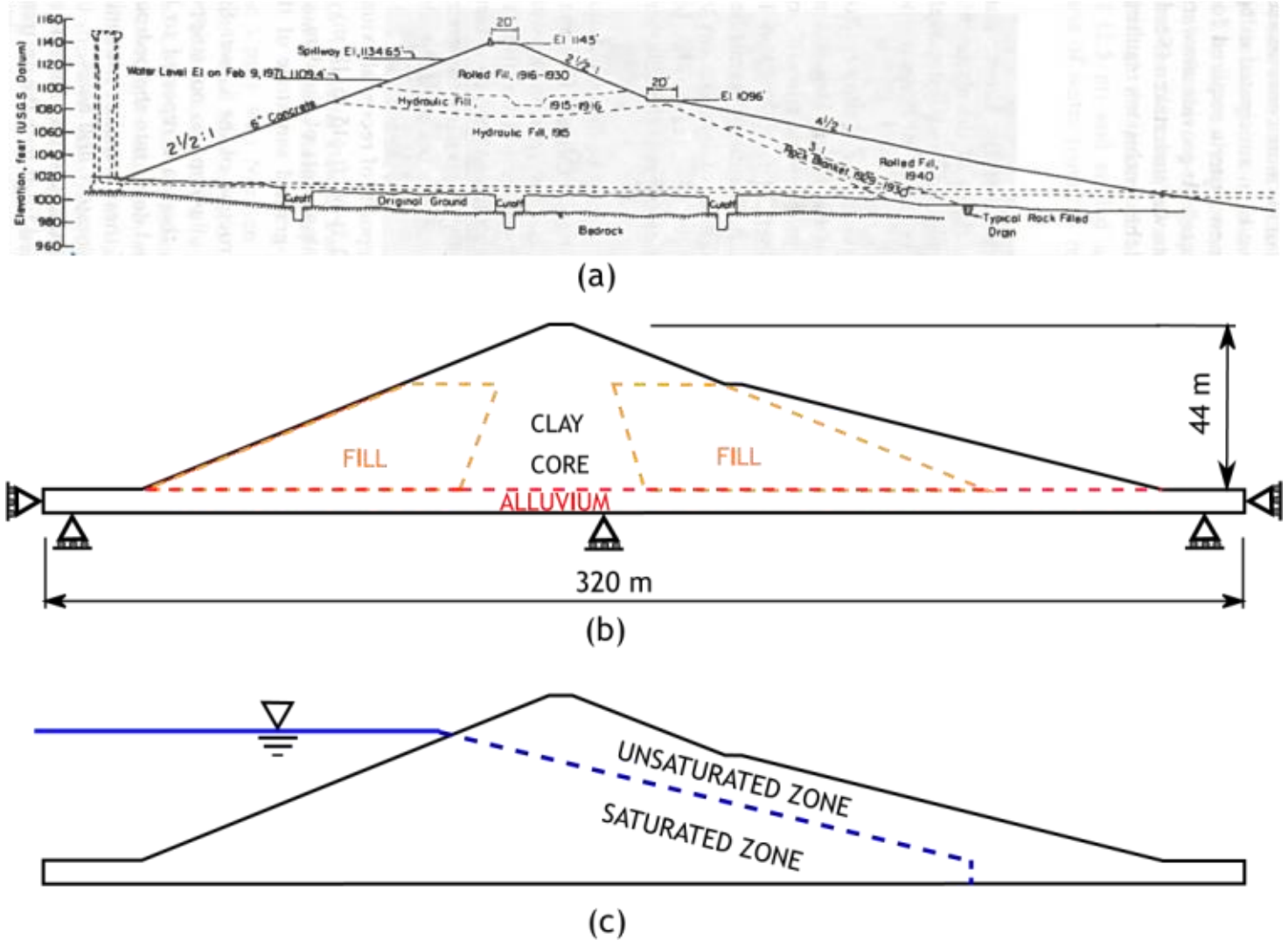

Figure 22: (a) Schematic of the dimensions of the dam reported in the literature, (b) Geometry, soil layers within the dam and boundary conditions in this study, and (c) sketch of water table in the reservoir and phreatic line in the dam.

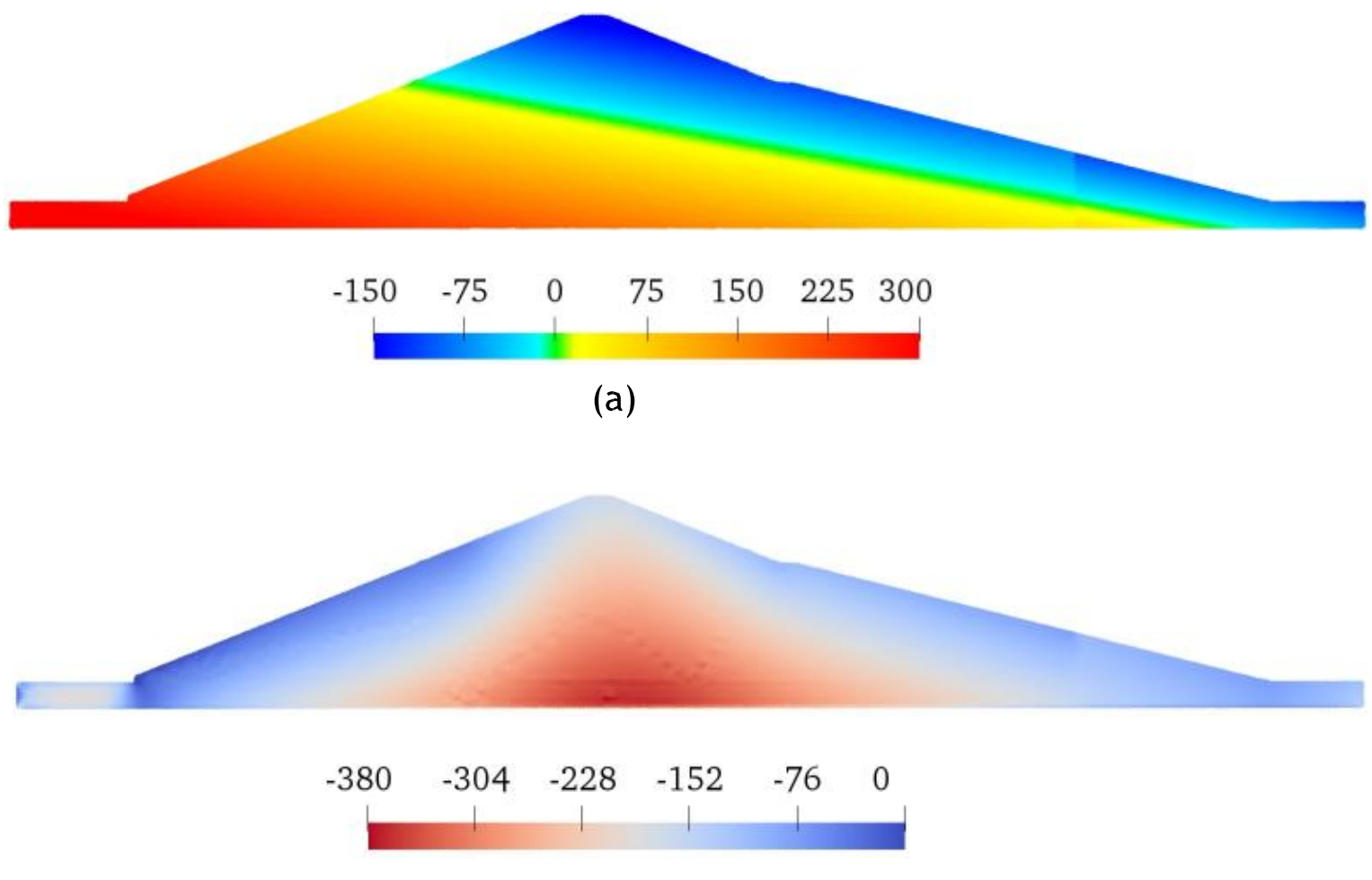

Figure 23: Contours of (a) the initial water pressure (kPa), and (b) the initial vertical effective stress (kPa) in the numerical model.

The bulk of the dam is made of various types of fines [12]. The core material is simplified as a single slightly overconsolidated clay in this example. The remainder of the dam is composed of looser and weaker hydraulic filler material outlined in brown in Figure 22. The input parameters for both materials are $\rho_s = 2 \times 10^3$ kg/m$^3$, $\rho_w = 1 \times 10^3$ kg/m$^3$, $\mu_w = 1 \times 10^{-3}$ Pa·s, $\varphi_0 = 0.375$ [1], $K = 2.33 \times 10^4$ kPa, $\mu_s = 1.167 \times 10^4$ kPa, $K_w = 2 \times 10^5$ kPa, $n = 3$, $s_a = 1.2 \times 10^3$ kPa. The plastic material parameters [32] for the core material are $p_c = -700$ kPa, $N = 2.0$, $M = 1.2$, $\kappa = 0.02$ and $\lambda = 0.09$. For the filler material the parameters are $p_c = -400$ kPa, $N = 1.75$, $M = 0.9$, $\kappa = 0.04$ and $\lambda = 0.13$.

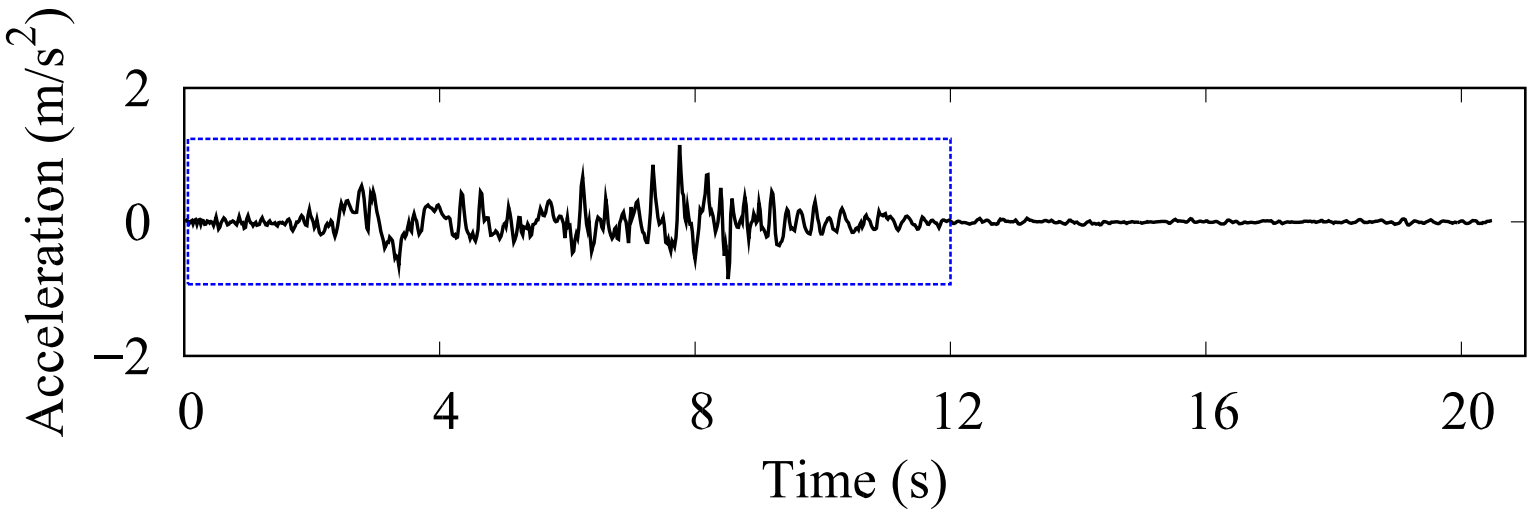

Figure 24: Recorded acceleration profile during the earthquake (Pacoima Dam accelerogram, www.strongmotioncenter.org).

Figure 23 (a) and (b) plots the initial geostatic stress and initial water pressure within the dam, which are generated by a quasi-static elastic analysis. With the crest of the dam taken as the datum, the initial effective skeleton stress is computed by the gravitational load on the soil within the dam. Pore water pressure along the upstream slope of the dam is prescribed by $p_w = \rho_w g h$, where h is the distance between the water table and the upstream slope. Pore water pressure and matric suction within the dam are computed through the seepage analysis as in [1]. The negative water pressure in the vadose zone is prescribed by $p_w = -\rho_w g h$, where $h$ is the vertical distance above the phreatic line.

Figure 24 plots the acceleration profile recorded during the earthquake. The strong ground motion only lasted for 12 s. The acceleration profile in the blue frame marks the input data in this study. In the coupled analysis, this acceleration profile is applied to the bottom of the dam as shown in Figure 22. The results are presented in Figures 25 and 26.

Figure 25 plots the snapshots of the equivalent plastic shear strain within the dam at (a) $t$ = 4 s, (b) $t$ = 10 s, and (c) $t$ = 16 s respectively. The results show that the large plastic shear strain has occurred at the toe of the upstream slope and the location beneath the core of the dam. The results may imply that failure initiated at the toe and progressed toward the core of the dam. Figure 26 presents the contours of water pressure within the dam at the same three time steps. The results in Figure 26 show that the water pressure increased at the base of the dam, which may imply the liquefaction of the soil under the dam core during the earthquake. The numerical results are consistent with the results based on the mixed finite element method in [1]. For a more realistic simulation of the dynamic failure of the dam under earthquake loads, a more representative constitutive model for different zones of soils within the dam will be required, which is beyond the scope of the present work.

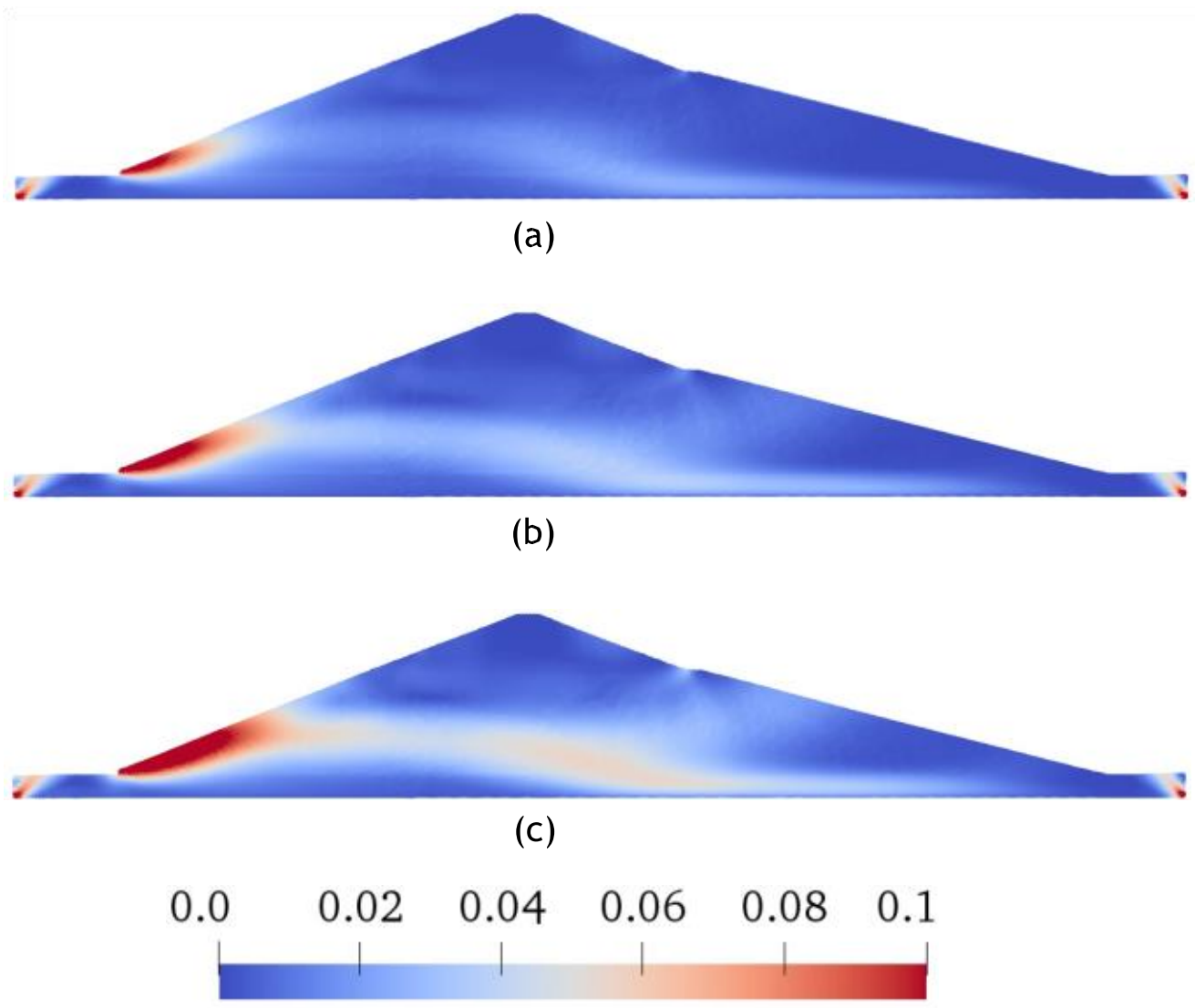

Figure 25: Snapshots of the contours of the plastic shear strain within the dam at (a) $t$ = 4 s, (b) $t$ = 10 s, and (c) $t$ = 16 s.

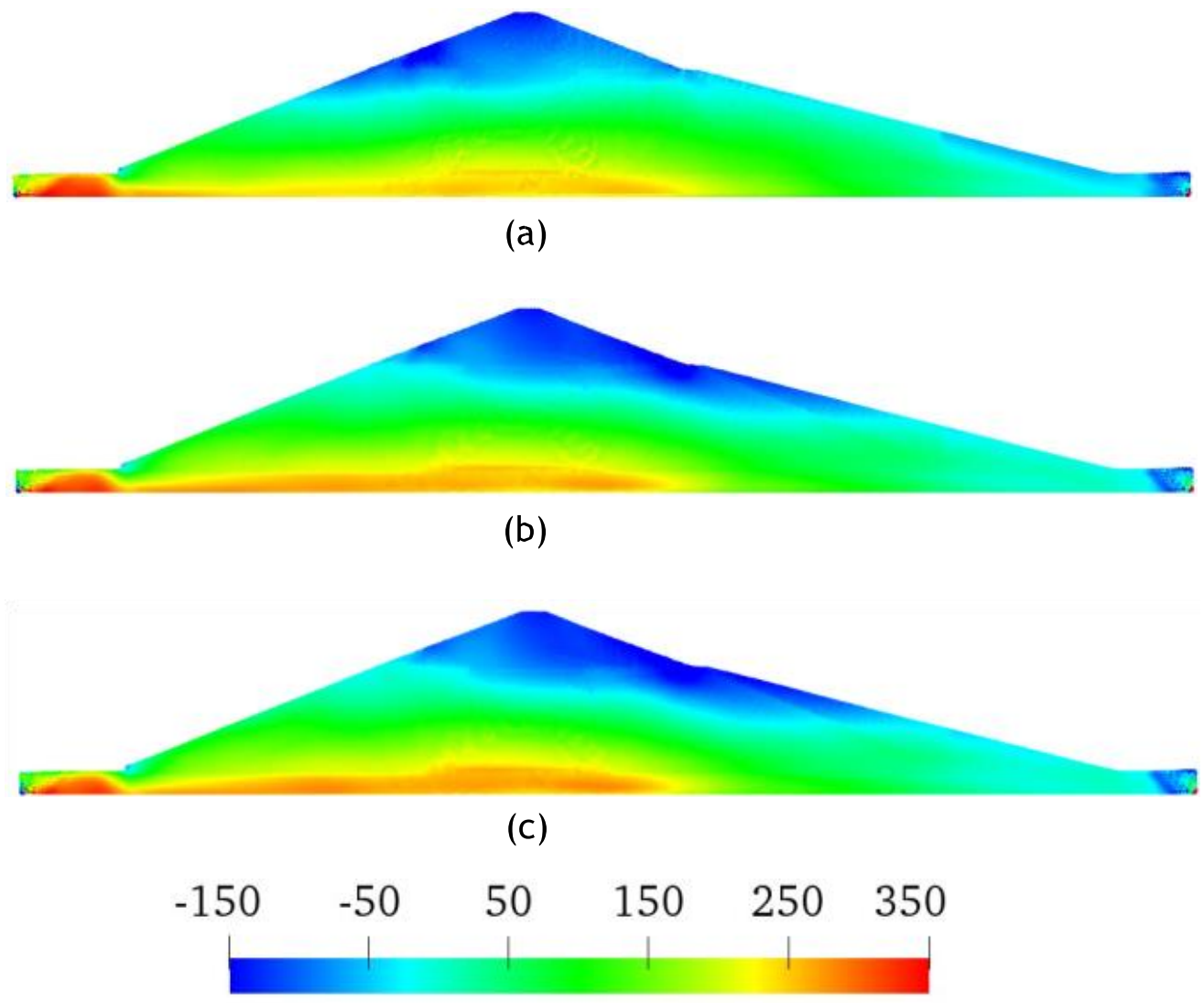

Figure 26: Snapshots of the contours of the water pressure (kPa) within the dam at (a) $t$ = 4 s, (b) $t$ = 10 s, and (c) $t$ = 16 s.

## 5. Closure

In this study, we have formulated a computational coupled large-deformation periporomechanics paradigm assuming passive air pressure for modeling dynamic failure and fracturing in variably saturated porous media. In this new computational periporomechanics paradigm, the coupled governing equations for both bulk and fracture material points are formulated in the current configuration through the updated Lagrangian-Eulerian framework. It is hypothesized that the horizon of a mixed material point remains spherical, and its neighbor points are determined in the current configuration. As a significant novelty, the mixed phreatic material points across the phreatic line are explicitly considered through the mixed peridynamic state concept. To incorporate the classical constitutive models into the nonlocal framework, we have developed the coupled constitutive correspondence principle with stabilization in the updated Lagrangian-Eulerian framework for the bulk and phreatic material points. We have numerically implemented the coupled large-deformation periporomechanics paradigm through an implicit-implicit fractional-step algorithm in time and a hybrid updated Lagrangian-Eulerian meshfree method in space. We first present numerical examples to validate the implemented coupled stabilized scheme and the fluid flow across the phreatic line in partially saturated porous media. We then conduct numerical examples to demonstrate the robustness and efficacy of the proposed computational framework in modeling fracturing and failure in partially saturated deformable porous media under static and dynamic loads.

**Acknowledgements**

This work has been supported by the US National Science Foundation under contract numbers 1659932 and 1944009. The support is gratefully acknowledged.